\magnification=\magstep1
\baselineskip=1.3\baselineskip

\newdimen\epsfxsize \baselineskip=18pt \overfullrule=0pt
\font\tenmsb=msbm10 \font\sevenmsb=msbm7 \font\fivemsb=msbm5
\newfam\msbfam
\textfont\msbfam=\tenmsb \scriptfont\msbfam=\sevenmsb
\scriptscriptfont\msbfam=\fivemsb
\def\Bbb#1{\fam\msbfam#1}

\def\H{{\Bbb H}}
\def\R{{\Bbb R}}

 \def\frac#1#2{{#1\over #2}}
 \def\ol{\overline}

 \def\dist{\hbox {\rm dist}}

 \def\df{{\mathop {\ =\ }\limits^{\rm{df}}}}
\def\bone{{\bf 1}}
\def\wt{\widetilde}
\def\wh{\widehat}
\def\n{{\bf n}}
\def\prt{\partial}

\def\cS{{\cal S}}

\def\cD{{\cal D}}
\def\cA{{\cal A}}
\def\cC{{\cal C}}
\def\C{{\bf C}}
\def\Om{{\Omega}}
\def\om{{\omega}}
\def\b{\partial}
\def\E{{\bf E}}
\def\P{{\bf P}}
\def\Q{{\bf Q}}
\def\s{{\bf s}}

\def\EE{{\cal E}}

\def\pf{\noindent{\bf Proof}. }

\input epsf.tex
\newdimen\epsfxsize
\newdimen\epsfysize

\def\qed{{\hfill $\square$ \bigskip}}
\def\sqr#1#2{{\vcenter{\vbox{\hrule height.#2pt
        \hbox{\vrule width.#2pt height#1pt \kern#1pt
           \vrule width.#2pt}
        \hrule height.#2pt}}}}
\def\square{\mathchoice\sqr56\sqr56\sqr{2.1}3\sqr{1.5}3}

\centerline{\bf ON THE ROBIN PROBLEM IN FRACTAL DOMAINS}

\footnote{$\empty$}{\rm Research partially supported by NSF grants
DMS-0244737 and DMS-0303310.}

\vskip0.2truein

\centerline{\bf Richard F.~Bass, Krzysztof Burdzy \rm and \bf
Zhen-Qing Chen}

\bigskip
\vskip0.6truein

{\narrower\narrower \noindent{\bf Abstract}. We study
the solution to the Robin boundary problem for the
Laplacian in a Euclidean domain. We present some families of
fractal domains where the infimum is greater than 0, and some other
families of domains were it is equal to 0. We also give a new
result on ``trap domains'' defined in [BCM], i.e., domains where
reflecting Brownian motion takes a long time to reach the center of
the domain.

 }

\bigskip
\noindent{\bf 1. Introduction}.

The Robin problem (also known as the ``third'' boundary problem)
for a Euclidean domain $D\subset \R^d$ is to find a function $u$
such that
$$\eqalignno{ \Delta u(x)&=0, \qquad x\in D,&(1.1)\cr
\frac{\partial u}{\partial {\bf n}}&=cu, \qquad x\in \partial
D,&(1.2)\cr}$$ with one or more side conditions, where $\n$ is the
unit inward normal vector field on $\partial D$, $\partial
u/\partial {\bf n}$ is the normal derivative of $u$ in the
distributional sense and $c>0$ is a constant. See Gustafson and
Abe [GA] for the history of this problem.

Our interest in the Robin problem stems from some recent
applications in physics, electrochemistry, heterogeneous catalysis
and physiology; see [FSF], [FS], [GFS], [Sa] and the references
therein. Consider the mixed Dirichlet-Robin problem
$$\eqalignno{ \Delta u(x)&=0, \qquad x\in D\setminus B_*,&(1.3)\cr
\frac{\partial u}{\partial {\bf n}}&=cu, \qquad x\in \partial D,&(1.4)\cr}$$
together with the side condition
$$ u(x)=1, \qquad x\in \prt B_*, \eqno (1.5)$$
where $B_*\subset D$ is a fixed closed ball with non-zero radius.
The solution to (1.3)-(1.5) represents the
steady state of a system in which some particles move randomly in
$D\setminus B_*$ and cross a semi-permeable membrane $\prt D$. The
other part of the boundary, $ \prt B_*$, is a source of particles
and can be controlled so that we can assume a condition of type
(1.5). The constant $c$ in (1.4) is a physical characteristic of
the membrane $\prt D$. One could consider a model with $c$
dependent on $x\in \prt D$ but we will not do that in the present
article. The constant $c$ will play no role in our theorems so
we will take $c=1$ in the rest of the article.

In some applied situations, it is desirable to have as much flux
through the boundary as possible. The points of a man-made or
natural membrane $\prt D$ where there is no flux can be considered
an inefficient use  of material. Hence, it is interesting to know
when the flux is non-negligible through all points of the
membrane. In other words, we would like to know whether $\inf
_{x\in \prt D} \prt u/\prt \n (x)>0$. In view of the relation
(1.4) between the flux $\prt u/\prt \n$ and the density $u$ of
particles and the maximum principle for the harmonic function $u$,
this condition is equivalent to
$\inf _{x\in D\setminus B_*} u(x) >0$.
(By Lemma 2.4 below, we know $u$ is non-negative.)

\bigskip
\noindent{\bf Definition 1.1}. {\sl We say that the whole surface
of $D$ is active if
  $$\inf _{x\in D\setminus B_*} u(x) >0.\eqno(1.6)$$
If it is not the case that the whole surface is active, we say
part of the surface is nearly inactive.
   }

\bigskip

In this paper we investigate the following problem.

\bigskip
\noindent{\bf Problem 1.2}. {\sl Give necessary and sufficient
conditions of a geometric nature for the whole surface of $D$ to be
active. }

\bigskip

It is not difficult to show that the whole surface of a bounded
Lipschitz domain is always active (see Remark 2.5(ii) below). We
have posed Problem 1.2 in terms of $u$ rather than $\prt u/\prt
\n$ because we are interested in non-Lipschitz domains $D$; so
there are some boundary points where $\n$ is not well-defined
while the solution $u$ is always well-defined, and, in fact, is
smooth in $D\setminus B_*$. We do not have a complete solution to
Problem 1.2, but we give a fairly explicit answer for some natural
families of domains with fractal boundary.

We will approach Problem 1.2 using probabilistic methods. This
agrees well with the motivating physical models. Suppose that $X$
is reflecting Brownian motion in $D$, $L$ is its local time on
$\prt D$, and $T_{B_*}$ is the hitting time of $B_*$ by $X$. When
$D$ is a bounded $C^3$-smooth domain, it is known that (see [MS]
and [Pa])
$$  u(x) = \E_x \left[ \exp \left( - {1\over 2} L_{T_{B_*}}\right) \right].
    \eqno(1.7)
$$
This formula indicates that the third boundary problem (1.4) is
more difficult to study from the probabilistic point of view than
the corresponding Dirichlet and Neumann problems. This is because
the Dirichlet problem corresponds to killed Brownian motion and
killing on the boundary presents no technical problems.  The
Neumann boundary problem corresponds to reflecting Brownian
motion. The construction of reflecting Brownian motion in an
arbitrary domain $D$ is a major technical challenge. Although this
feat has been accomplished long time ago by Fukushima [Fu] on an
abstract compactification, called the Martin-Kumarochi
compactification, of $D$, many questions about the construction of
reflecting Brownian motion on the Euclidean closure of a domain
remain open (see [BBC]). Formula (1.7) shows that the Robin
boundary problem (1.3)-(1.5) requires the construction and
understanding of the local time. This is harder than constructing
reflecting Brownian motion itself, because it is known that
reflecting Brownian motion does not have a semimartingale
decomposition in some domains. For some results in this area, see,
e.g., DeBlassie and Toby [DT]. For information on the eigenvalue
problem for the Laplacian with Robin boundary conditions, see
Smits [Sm1], [Sm2].

\bigskip

The following are the main contributions of this paper. The list
includes some technical results that may have independent
interest.

\item{(i)} The solution of Problem 1.2 for a class of domains with
fractal boundaries (Theorems 3.2 and 4.3).

\item{(ii)} A characterization of a class of ``trap domains'' in
dimensions 3 and higher, improving a result in [BCM] (Theorem
5.1).

\item{(iii)} Clarification of the rigorous meaning of
solution to the differential equation (1.3)-(1.5), its existence,
uniqueness, and probabilistic representation for
non-smooth domains ((2.3) and Lemma 2.4).
In particular, we show that the solution to (1.3)-(1.5) is non-negative.

\item{(iv)} A semimartingale decomposition of reflecting Brownian
motion in a class of fractal domains (Theorem 2.2).

\item{(v)} A sharp estimate for the Green function with Neumann
boundary conditions in long and thin domains (Lemma 4.4).

\item{(vi)} A new version of the Neumann boundary Harnack
principle, stronger than the one in [BH] (Lemma 2.8).

\item{(vii)} The proof that reflecting Brownian motion starting
from the cusp point is not a semimartingale, for some cusps
(Remark 4.14). This complements a result of Fukushima and Tomisaki
[FT].

\bigskip

A simple example illustrating our main theorems is a cusp domain,
defined for a fixed $\alpha >1$ by
$$
 D= \left\{x = (x_1, x_2, \dots, x_d): 0 < x_1 < 1
\hbox{ and } x_1^\alpha >(x_2^2 + \dots + x_d^2)^{1/2}\right\}.
 $$
 Applying the main results (Theorems 3.2 and  4.3) of this paper, we
 show in Example 3.4 (for $d=2$) and Example 4.13 (for $d\geq 3$)
that the whole boundary of $D$ is active if $\alpha \in (1, 2)$,
and part of $\partial D$ is nearly inactive if $\alpha \geq 2$.
There are more examples given in Sections 3 and 4.
\bigskip

The paper is organized as follows. Section 2 contains some
technical preliminaries, many of which may have independent
interest. Section 3 presents the solution to Problem 1.2 for a
class of 2-dimensional domains, using techniques developed in
[BCM]. Section 4 is devoted to Problem 1.2 in dimensions 3 and
higher. Finally, Section 5 presents an application of the
techniques developed in Section 4 to ``trap domains'' in
dimensions 3 and higher.

\bigskip

We would like to thank Marcel Filoche, Masatoshi Fukushima, Don
Marshall, Stephen Rohde, Bernard Sapoval and Tatiana Toro for
stimulating discussions and valuable advice.

\bigskip

\noindent {\bf 2. Reflecting Brownian motion in domains with
fractal boundaries and the Neumann boundary Harnack principle}.

This section is devoted to two important technical aspects of this
paper. First, we will show that reflecting Brownian motion has a
semimartingale decomposition for a class of fractal domains that
contains some natural examples. The second technical result is a
boundary Harnack principle for harmonic functions satisfying
Neumann boundary conditions.

We will let $|\,\cdot\,|$ stand for the Euclidean norm in $\R^d$
(for any dimension $d\geq 1$), for the volume ($d$-dimensional
Lebesgue measure) of a set $A\subset \R^d$, and for the
$(d-1)$-dimensional surface area of the boundary $\prt A$ of a set
$A\subset \R^d$. The meaning will be obvious from the context so
this notation should not lead to any confusion. For an open set
$D$ of $\R^d$, $C_c(D)$ and $C_c^\infty (D)$ denote the space of
continuous functions with compact support in $D$ and the space of
smooth functions with compact support in $D$, respectively.

A ball with center $x$ and radius $r$ will be denoted $B(x,r)$.
The notation will refer to an open ball, unless noted otherwise.

The harmonic measure of a set $A\subset \b D$ in the domain $D$,
relative to $z$, will be denoted $\om(z,A,D)$.

The distribution of Brownian motion in $D\setminus B_*$ starting
from $x \in \ol{D\setminus B_*}$, reflected on $\prt D$, and
killed on $\prt B_*$ will be denoted $\P_x$. The corresponding
expectation will be denoted $\E_x$. The hitting time of a set $A$
will be denoted $T_A$, i.e., $T_A = \inf\{t>0: X_t\in A\}$. We
will sometimes write $T^X_A$ or $T^Y_A$ to show the dependence of
the hitting time on the process.

We will use elements of excursion theory and Doob's $h$-processes.
See [D] for the discussion of $h$-transforms in the case of
(non-reflecting) Brownian motion, and [Sh] for conditioning of
general Markov processes. Elements of excursion theory can be
found in [Mv], [Bl], [Bu] and [Sh].

A real-valued function $f$ defined on $A\subset \R^d$ is called
Lipschitz with constant $\lambda< \infty$ if $|f(x)-f(y)|\leq
\lambda |x-y|$ for all $x,y\in A$. A domain $D$ is called
Lipschitz if there exist $r>0$ and $\lambda< \infty$ such that for
every $x\in \prt D$, the set $\prt D \cap B(x,r)$ is the graph of
a Lipschitz function with constant $\lambda$ in some orthonormal
coordinate system. We call $(\lambda,r)$ the Lipschitz
characteristics of $D$.

\bigskip
\noindent{\bf Definition 2.1}. We will say that a domain $D$
belongs to class $\cD$ if there exists an increasing sequence of
domains $D_n\subset D$ with the following properties.

\item{(i)} Each $D_n$ is a Lipschitz domain with characteristics
$(\lambda, r_n)$ (all the $\lambda$'s are the same but the $r_n$'s
may differ) and $\bigcup_{n=1}^\infty D_n=D$.

\item{(ii)} For every $n\geq 1$, the set $\prt D_n \cap \prt D$ is
a subset of the relative interior of $\prt D_{n+1} \cap \prt D$.

\item{(iii)}  $\sup_{n\geq 1} |\prt D_n| < \infty$ and
$\lim_{n\to \infty} |\prt D_n \setminus \prt D| =0$.

The set $\prt_L D \df \bigcup_n \prt D_n \cap \prt D$ will be
called the Lipschitz part of $\prt D$.

\bigskip

Every bounded Lipschitz domain is in $\cD$. See Examples 3.4, 3.6,
4.13 and 4.14 below for domains $D\in \cD$ which are not Lipschitz.

Constructing a reflecting Brownian motion on a non-smooth domain
$D$ is a delicate problem. Let
$$ W^{1,2}(D)\df \{f \in L^2(D, dx): \, \nabla f \in L^2(D, dx)\}
$$
be the Sobolev space on $D$ of order $(1, 2)$. Fukushima [Fu] used
the Martin-Kuramochi compactification $D^*$ of $D$ to construct a
continuous diffusion process $X^*$ on $D^*$ with transition
semigroup denoted $P_t$, such that
$$ \{ f \in L^2(D, dx): \sup_{t>0} {1\over t} \int_D f(x)
(f(x)-P_t f(x)) dx<\infty\}=W^{1,2}(D)
$$
and for $f\in W^{1, 2}(D)$,
$$ \EE(f, f)\df \lim_{t\to 0} {1\over t} \int_D f(x)
(f(x)-P_t f(x)) dx = {1\over 2} \int_D |\nabla f (x)|^2 dx.
$$
The pair $(\EE, W^{1,2}(D))$ is called the Dirichlet space of
$X^*$ in $L^2(D^*, m)$, where $m$ is Lebesgue measure on $D$
extended to $D^*$ by setting $m(D^*\setminus D)=0$. See [FOT] for
definitions and properties of Dirichlet spaces, including the
notions of quasi-everywhere, quasi-continuous, etc. The process
$X^*$ could be called reflecting Brownian motion in $D$ but it
lives on an abstract space $D^*$ that contains $D$ as a dense open
set. Chen [C1] proposed referring to  the quasi-continuous
projection $X$ of $X^*$ from $D^*$ into the Euclidean closure
$\overline D$ as reflecting Brownian motion in $D$. The projection
process $X$ is a continuous process on $\overline D$, but in
general $X$ is not a strong Markov process on $\overline D$ (for
example this is the case when $D$ is the unit disk with a slit
removed). However when $D$ is a Lipschitz domain, it is shown that
$X$ is the usual reflecting Brownian motion in $D$ as constructed
in [BH].
It was proved in [C1] that, roughly speaking, if $\partial D$ has
``finite surface measure,'' then $X$ is a semimartingale and has a
Skorokhod decomposition.
This result was further sharpened in [CFW]. See the introductions
 of [C1] and [CFW] for  the history of constructing reflecting
Brownian motion on non-smooth domains.

\bigskip
\noindent{\bf Theorem 2.2}. {\sl If $D\in \cD$, then reflecting
Brownian motion $X$ in $D$ starting from $x\in D\cup \prt_L D$ has
a semimartingale decomposition $X_t = x + W_t + N_t$, where $W_t$
is a $d$-dimensional Brownian motion,
 $$N_t = \int_0^t \n(X_s) dL_s,
 $$
and $L$, the local time, is a non-decreasing continuous process
that does not increase when $X$ is not in $\prt _L D$, i.e.,
$\int_0^\infty \bone_{(\prt_L D)^c} (X_t) dL_t = 0$. The Revuz
measure of $L$ for the process $X^*$ is surface measure on
 $\prt_L D$.}

\bigskip

Note that the local time $L$ in our theorem satisfies the
condition $\int_0^\infty \bone_{(\prt_L D)^c} (X_t) dL_t = 0$,
which is stronger than the usual condition $\int_0^\infty \bone_D
(X_t) dL_t = 0$.

\bigskip
\noindent {\bf Proof}. Let $\{D_n, n\geq 1\}$ be the increasing sequence
of Lipschitz domains in the definition of $D\in {\cal D}$.
Let $D^*$ be the Martin-Kuramochi compactification
of $D$ used in [Fu]. To be precise, for every $\alpha >0$, let
${\cal H}_\alpha$ denote the space of all $h$
in $D$ such that $(\alpha - {1\over 2} \Delta ) h=0$ in $D$
and having
$$ \EE_\alpha (h, h) \df {1 \over 2} \int_D |\Delta h (x) |^2 \, dx
 + \int_D u(x)^2 \, dx < \infty.
$$
For $y\in D$, let $x\mapsto H_\alpha (x, y)$ be the unique
$\alpha$-harmonic function in ${\cal H}_\alpha$ such that
$$ \EE_\alpha (H_\alpha (\, \cdot\, ,y ), \, v(\,\cdot\,)) = v(y)
\qquad \hbox{for every } v \in {\cal H}_\alpha.
$$
Let $G^0_\alpha (x, y)$ be the $\alpha$-resolvent density
function for Brownian motion in $D$ killed upon exiting $D$.
Define
$$ G_\alpha (x, y)\df G^0_\alpha (x, y) + H_\alpha (x, y).
$$
It is shown in [Fu] that $x\mapsto G_\alpha (x, y)$ is continuous
on $D\setminus \{y\}$ and $G_\alpha (x, y)=G_\alpha (y, x)$.
Define a metric $\delta$ on $D$ by
$$ \delta (x, y) = \int_D
 (| G_1 (x, z)-G_1 (y, z)| \wedge 1) \, dz
$$
and let $D^*$ be the completion of $D$ under the metric $\delta$.
Fukushima [Fu] showed that there is a conservative continuous Hunt
process $X^*$ on $D^*\setminus N$ associated with the Dirichlet
space $(\EE, W^{1,2}(D))$ on $L^2(D^*, m)$, where $N$ is a set
that has zero capacity with respect to $(\EE, W^{1,2}(D))$ and $m$
is Lebesgue measure on $D$ extended to $D^*$ by defining
$m(D^*\setminus D)=0$. Since each coordinate function $x_i \in
W^{1, 2}(D)$, then each coordinate function
 admits a quasi-continuous version on $D^*$, which
will be denoted as $f_i$. Note that $(f_1, \cdots, f_d)$ is
defined quasi-everywhere on $D^*$ and is a quasi-continuous map
from  $D^*$ into $\overline D$. Define
$$ X=(f_1(X^*), \cdots, f_d(X^*)).
$$
Then $X$ is a conservative continuous process on $\overline D$,
which is called reflecting Brownian motion on $\overline D$ in
[C1]. It coincides with the usual reflecting Brownian motion when
$D$ is a bounded Lipschitz domain.

Let $X^n$ be reflecting Brownian motion on $\overline D_n$. It is known
from [BH] that $X^n$ has a H\"older continuous transition density function
$p^n(t, x, y)$ on $(0, \infty) \times \overline D_n \times \overline D_n$.
Its $\alpha$-resolvent density function will be denoted as $G^n(x, y)$.
Define
$$ \tau_n = \inf\{ t\geq 0: X^n_t \in \partial D_n \setminus \partial D \}
$$
and
$$ G_\alpha^{n, 0} (x, y_0)= G_\alpha (x, y_0)-\E_x^{n}
\left[ e^{-\alpha \tau_n} G_\alpha (X^n_{\tau_n}, y_0) \right].
$$
It is easy to verify that $G_\alpha^{n , 0}$ is the $\alpha$-resolvent
for reflecting Brownian motion in $\overline D_n$ killed upon
hitting $\partial D_n \setminus \partial D$. Thus we have
$$ G_\alpha (x, y_0)=G^n_\alpha (x, y)+\E_x^{n}
\left[ e^{-\alpha \tau_n} \left( G_\alpha (X^n_{\tau_n}, y_0)
 - G_\alpha^n (X^n_{\tau_n}, y_0) \right) \right].
$$
By [BH], $x\mapsto G^n_\alpha (x, y_0)$ is continuous on
$\overline D_n \setminus \{y_0\}$ and
$$ x\mapsto \E_x^n \left[ e^{-\alpha \tau_n} \left( G_\alpha (X^n_{\tau_n}, y_0)
 - G_\alpha^n (X^n_{\tau_n}, y_0) \right) \right]
$$
is continuous on $\partial D_n\cap  \partial D$ since it is
harmonic in $D_n\setminus \{ y_0 \}$ with zero Neumann boundary
conditions on $\partial D_n\cap
\partial D$ and zero Dirichlet boundary conditions on $\partial D_n \setminus
\partial D$. Hence we conclude that $x\mapsto G(x, y_0)$ extends continuously
to $\partial D_n \cap \partial D$ under the Euclidean topology for
every $n\geq 1$ and hence to $\partial_L D=\bigcup_{n=1}^\infty
(\partial D_n \cap \partial D)$. This implies that
$$\partial_L D\subset (D^*\setminus D) \cap \partial D.
$$
Note that on $\partial_L D$, $(f_1, \cdots, f_d)$ is the identity map
and so $X_t^*=X_t$ when $X^*_t\in \partial_L D$.

Let $\sigma_k$ denote surface measure on $\partial D_k$ and let
${\bf n}_k(x)$ be the unit inward normal vector field on $\partial
D_k$ which is defined almost everywhere with respect to
$\sigma_k$. By the definition of $D\in {\cal D}$,
$$ k \mapsto \sigma_k (\partial D_k \cap \partial D )
$$
is an increasing function  and
$$ \lim_{k\to \infty} \sigma_k (\partial D_k) =
   \lim_{k\to \infty} \sigma_k (\partial D_k \cap \partial D)
  =\sigma ( \partial_L D),
$$
since  $\partial D_k \cap \partial D \subset \partial D_{k+1} \cap
\partial D$ and $\lim_{k\to \infty} \sigma_k (\partial D_k
\setminus \partial D)=0$. Here $\sigma$ is surface measure on
$\partial_L D$. Since $\sup_{k\geq 1} \sigma_k (\partial D_k) <
\infty$, there exist a subsequence $\{k_j, j\geq 1\}$ and finite
signed measures $(\nu_1, \cdots , \nu_d)$ on $D^*$ such that
$\n_{k_j}\sigma_{k_j}$ converges weakly on $D^*$ to $(\nu_1,
\cdots , \nu_d)$; that is,
$$ \lim_{j\to \infty} \int_{D^*} (g_1(x), \cdots , g_d(x))
 \cdot \n_{k_j}(x)\, \sigma_{k_j}(dx)
= \sum_{i=1}^d \int_{D^*} g_i(x)\, \nu_i(dx), \eqno(2.1)
$$
for all bounded continuous functions $\{g_1, \cdots, g_d\}$ on $D^*$.
For every $1\leq i \leq d$ and $k\geq 1$,
$$\eqalign{
 |\nu_i|( D^* \setminus (\partial D_k \cap \partial D))
& \leq \lim_{j\to \infty} \sigma_{k_j}
  ( D^* \setminus (\partial D_k \cap \partial D))  \cr
&  = \lim_{j\to \infty} \sigma_{k_j}
  ( \partial D_{k_j}  \setminus (\partial D_k\cap \partial D)) \cr
& = \sigma (\partial_L D\setminus \partial D_k) .\cr}
$$
Thus for $1\leq i \leq d$,
$$ |\nu_i|(D^*\setminus \partial_L D)
 = \lim_{k\to \infty} |\nu_i|( D^* \setminus (\partial D_k \cap \partial D))
\leq \lim_{k\to \infty} \sigma (\partial_L D \setminus \partial D_k)=0.
\eqno(2.2)
$$
On the other hand, by the definition of $D\in {\cal D}$,
${\bf n}_k \,\sigma_k$ converges weakly on $ \overline D$ to
${\bf n}\, \sigma$,
where  ${\bf n}$ is the unit inward normal vector
field of $D$ on $\partial_L D$ in the following sense:
$$ \lim_{k\to \infty} \int_{\overline D} (g_1(x), \cdots , g_d(x))
 \cdot \n_k(x)\, \sigma_k(dx) = \int_{\overline D} (g_1(x), \cdots , g_d(x))
 \cdot \n(x) \,\sigma(dx)
$$
for all bounded continuous functions $\{g_1, \cdots, g_d\}$ on $\overline D$
that vanish
on $\partial D_n \setminus \partial D$ for some $n\geq 1$.

Since $\partial_L D \subset D^*\cap \overline D$, we conclude
from (2.1) and (2.2) that
$$ (\nu_1, \cdots , \nu_d)= \n \, \sigma \qquad \hbox{on } D^*.
$$

By  Theorem 4.4 of [C1], $\sigma$ is a smooth measure of $X^*$ and
thus it determines a positive continuous additive function $L$ of
$X^*$. Moreover,
$$ X_t=X_0+W_t+\int_0^t \n (X_s) dL_s \qquad \hbox{for } t\geq 0,
$$
where $W$ is a $d$-dimensional Brownian motion. The above
Skorokhod decomposition holds for quasi-every starting point
$X^*_0$ in $D^*$ with $X_0=f(X_0)$. However, since the
$\alpha$-resolvent density function $x\mapsto G(x, y_0)$ is
continuous on $\partial_L D \cup (D\setminus \{y_0\})$, reflecting
Brownian motion $X^*$ can be defined to start from every point
$x\in D \cap \partial_L D$ (cf. [FOT]). Hence the above Skorokhod
decomposition holds for every starting point $X_0 \in D \cup
\partial_L D$. Clearly, since $\sigma$ is carried on $\partial_L
D$,
$$ \int_0^\infty 1_{\{ X_s\notin \partial_L D\}} dL_s
  = \int_0^\infty 1_{\{ X^*_s\notin \partial_L D\}} dL_s
  =0.
$$
This proves the theorem.
\qed

\bigskip
\noindent{\bf Remark 2.3}. Let $\tau_{D\cup \prt_L D}$ be the
first exit time from $D\cup \prt_L D$ by reflecting Brownian
motion on $D$. Starting from $x \in D\cup \partial_L D$, $\{X_t,
t<\tau_{D\cup \partial_L D}\}$ is a strong Markov process on
$D\cup \partial_L D$, since it coincides with $\{X^*_t,
t<\tau_{D\cup \partial_L D}\}$. Here
$$\tau_{D\cup \partial_L D} \df \inf\{ t>0: X_t \notin D\cup \partial_L D\}
= \inf\{ t>0: X^*_t \notin D\cup \partial_L D\} .
$$
However even under the conditions of Theorem 2.2, reflecting
Brownian motion on $D$ may not be a strong Markov process. For
example, let $D$ be the union of $\{(x, y)\in \R^2: |y|>|x| \hbox{
and } |y|\leq 1\}$ and $\{ (x, y)\in \R^2: 1<x^2+y^2<4\}$. Then
clearly reflecting Brownian motion $X$ on $\overline D$ can not have the strong
Markov property since when $X_t$ is at the origin {\bf 0}, one can
not tell how it will be reflected unless one knows where it came
from. Of course, in this example, for starting points in
$\overline D$ other than the origin ${\bf 0}$, reflecting Brownian
motion will not visit ${\bf 0}$. But one can modify this example
so that the set of such non-Markovian points has positive capacity
so it will be visited by the reflecting Brownian motion. Here is
such an example. Let $K$ be the standard Cantor set in $[0, 1]$.
Let $C=\{(x, y)\in \R^2: \, |x|<|y|<1\}$. Define $D\subset \R^2$
by
$$ D=  \{(x, y): 1<x^2+y^2<9\}
\cup \bigcup_{x\in K} \left( (x, 0)+C) \right) .
$$
Clearly $D$ satisfies the assumptions of Theorem 2.2. Note $\wt
K\df K\times \{0\}$ is the set of non-Markovian points and $\wt K$
has positive capacity (see [C2]) and so will be visited by
reflecting Brownian motion in $D$. \qed

\bigskip

Now we make precise the meaning of solution to the partial
differential equation with Robin and Dirichlet boundary conditions
(1.3)-(1.5) for $D\in \cD$. Define
$$ W^{1,2}(D; B_*)\df \{u\in W^{1,2}(D): \, u=0
\hbox{ q.e. on the closed ball } B_* \}.
$$
We say $u$ is a (weak) solution of (1.3)-(1.5) if the following
two conditions are satisfied.

\item{(i)} $u$ and its distributional derivative $\nabla u$ are in
$L^2(D\setminus B_*)$ and for any bounded $g\in W^{1,2}(D;
B_*)$,
$$   \int_{D\setminus B_*} \nabla g(x) \cdot \nabla u(x) \, dx
=- c \, \int_{\partial_L D} g(x) u(x) \sigma (dx). \eqno(2.3)
$$
\item{(ii)} $u$ is continuous in a neighborhood of $\partial B_*$ and
$u=1$ on $\partial B_*$.

\medskip
Note that any $f\in W^{1,2}(D; B_*)$ admits a quasi-continuous version
on $D^*\setminus B_*$. Throughout this paper, we will always represent such $f$
by its quasi-continuous version,
which will still be denoted as $f$. In particular, $f$ is well defined q.e. on
$D^*\setminus D$. Since $\partial_L D\subset D^*\setminus D$ and $\sigma$ is a
smooth measure of $X^*$ according to Theorem 2.2, $f$ is well defined
$\sigma$-a.e. on $\partial_L D$ for every $f\in W^{1,2}(D; B_*)$.
Hence the right hand side of (2.3) is well defined.

\bigskip

\noindent{\bf Lemma 2.4}. {\sl The partial differential equation
with Robin and Dirichlet boundary conditions (1.3)-(1.5),
where $c>0$ is a constant,
has a unique solution $u(x)$ given by $u(x) =
\E_x  \left[ \exp (-{c\over 2} L_{T_{B_*}}) \right]$.
In particular, $u$ is non-negative. }

\bigskip
\noindent {\bf Proof}.
We first establish existence.
Note that $u(x)\df \E_x [\exp(-{c\over 2}
L_{T_{B_*}})]$ is well defined for every $x\in D\cup \partial_L D$
and for q.e. $x\in D^*$. By the Markov property of $X^*$, for
$x\in D \cup \partial_L D$ (as well as for q.e. other $x\in D^*$),
$$
  v(x)\df 1-\E_x\left[ \exp (-{c\over 2} L_{T_{B_*}})\right] =
{c\over 2} \E_x \left[ \int_0^{T_{B_*}} u(X^*_s) dL_s \right] .
$$
Let $X^{*, 0}$  be reflecting
Brownian motion $X^*$ killed upon hitting $B_*$ and let the
transition semigroup be denoted by $\{P^0_t, t\geq 0\}$.
It is known (cf. [FOT]) that the Dirichlet form of $X^{*, 0}$ is
$(\EE, W^{1,2}(D; B_*))$ on $L^2(D^*\setminus B_*, m)$.
For q.e. $x\in D^*\setminus B_*$,
$$ v(x)-\P^0_tv(x)= {c\over 2}
\E_x \left[ \int_0^t u(X^{*, 0}_s) dL_s \right].
$$
Hence
$$\eqalign{
 \lim_{t\to 0} {1\over t} \int_{D\setminus B_*} v(x)
(v(x)-P^0_t v(x)) \, dx
& = {c\over 2} \lim_{t\to 0}   \int_{D\setminus B_*} v(x)
\E_x \left[  \int_0^t u(X^{*, 0}_s) dL_s \right] dx \cr
& ={c\over 2} \int_{\partial_L D} v(x) u(x) \sigma (dx) <\infty. \cr}
$$
Thus $v\in W^{1,2}(D; B_*)$, and a similar calculation to the above
yields that for any bounded $g\in W^{1,2}(D; B_*)$
$$\eqalign{
 {1\over 2} \int_{D\setminus B_*} \nabla g(x) \cdot \nabla v(x) \, dx
&= \EE(g, v)  \cr
&= \lim_{t \to 0} {1\over t} \int_{D\setminus B_*} g(x)
(v(x)-P^0_t v(x)) \, dx \cr
&= {c\over 2} \int_{\partial_L D} g(x) u(x) \sigma (dx). \cr}
$$
This shows that $v$ is harmonic in $D\setminus B_*$ and
${\partial v \over \partial \n} = -c u $.
In particular, $v$ is continuous in $D\setminus B_*$.
 Since every point of $\partial B_*$  is regular, we see that
$v$ vanishes continuously on $\partial B_*$. Translating these
properties to the function $u=1-v$
shows that $u$ is a solution to (1.3)-(1.5).

Now we show the uniqueness. Suppose that $u_1$ and $u_2$ are
two solutions for (1.3)-(1.5). Define $w\df u_1-u_2$.
Then $w\in W^{1,2}(D; B_*)$ and it follows from (2.3) that
$$   \int_{D\setminus B_*} \nabla g(x) \cdot \nabla w(x) \, dx
=- c \, \int_{\partial_L D} g(x) w(x) \sigma (dx).
$$
Letting $g=((-n)\vee w)\wedge n$ in the above and then letting $n\to \infty$,
we have
$$   \int_{D\setminus B_*} |\nabla w(x)|^2 \, dx
=- c \, \int_{\partial_L D} | w(x)|^2 \sigma (dx).
$$
Since $c>0$, we must have
$$ \int_{D\setminus B_*} |\nabla w(x)|^2 \, dx
= \int_{\partial_L D} | w(x)|^2 \sigma (dx) =0. \eqno (2.4)
$$
Since $D\setminus B_*$ is connected, $w$ has to be constant
in $D\setminus B_*$, while the second equality in (2.4)
 implies that
$w=0$ $\sigma$-a.e. on $\partial_L D$.
Therefore $w=0$ in $D\setminus B_*$ and hence $u_1=u_2$.
This establishes the uniqueness and completes the proof of this Lemma.
\qed

\bigskip

\noindent{\bf Remark 2.5}. (i) A simple modification of the above
argument establishes the existence and uniqueness for solutions to
the Robin problem (1.3)-(1.5) with $c$ being a bounded
non-negative function. The solution in this case can be
represented as
$$ u(x)=\E_x \left[ \exp \left(-{1\over 2} \int_0^{T_{B_*}} c(X_s) dL_s \right)
 \right]  \qquad \hbox{for } x\in D\setminus B_*.
$$

\medskip

(ii) Suppose that $D$ is a bounded Lipschitz domain in $\R^d$ with
$d\geq 3$ and $u$ is the solution to (1.3)-(1.5). By Jensen's
inequality, we have
$$ u(x) \geq \exp \left( -\frac{c}2 \E_x L_{T_{B_*}} \right)
\qquad \hbox{for } x\in D\setminus B_*.
$$
Let $G_{D\setminus B_*}$ be the Green function of the reflecting
Brownian motion killed upon hitting $B_*$. It is known from [BH]
that
$$ G_{D\setminus B_*} (x, y) \leq \frac{c_1}{|x-y|^{d-2}}
\qquad \hbox{for } x, y \in \overline D \setminus B_*.
$$
It follows then
$$ \sup_{x\in D\setminus B_*} \E_x
L_{T_{B_*}}
 =\sup_{x\in D\setminus B_*} \int_{\partial D} G_{D\setminus B_*}
   (x, y) \sigma (dy)
 \leq \sup_{x\in D\setminus B_*} \int_{\partial D}
   \frac{c_1}{|x-y|^{d-2}} \sigma (dy ) <\infty.
$$
Hence $\inf_{x\in D\setminus B_*} u(x) >0$. In other words, the
whole surface of a bounded Lipschitz domain in $\R^d$ with
$d\geq 3$ is always active. \qed

\bigskip

Let $D\subset \R^d$ be a Lipschitz domain and let $O$ be a
connected open set in $\R^d$. The following definition of
``Neumann boundary conditions'' for a harmonic function is
standard in analysis and PDE (cf. [K]).

\medskip

\noindent{\bf Definition 2.6}. A function  $h$ defined on $D\cap
O$ is said to be harmonic in $D\cap O $ with zero Neumann boundary
conditions on $ \prt D\cap O $ if $h\in W^{1,2}(O_1\cap D)$ for
every relatively compact open subset $O_1$ of $O$ and
$$ \int_{O\cap D} \nabla h (x) \cdot \nabla \psi (x) dx
  =0  \eqno(2.5)
$$
for every $ \psi \in C^\infty_c (O_1)$ and consequently
for every continuous $\psi \in W^{1,2} (O_1)$ that vanishes
on $\partial O_1$.

\medskip
The following lemma says that functions expressed in terms of the hitting distribution
of reflecting Brownian motion in $D$ are harmonic functions with zero Neumann boundary
conditions in the sense of Definition 2.6.

\bigskip

\noindent{\bf Lemma 2.7}. {\sl Let $X$ be reflecting Brownian
motion in the Lipschitz domain $D$, and $O$ a connected open
subset of $\R^d$. Define $\tau\df \inf\{t>0: \, X_t \notin
\overline D \cap O\}$. Then for any bounded measurable function
$\psi$ on $\partial O \cap \overline D$,
$$ h(x)\df \E_x \left[ \psi (X_{\tau}) \right], \qquad x \in \overline D \cap O,
$$
is a harmonic function in $D\cap O$ with zero Neumann boundary conditions
at $\partial D \cap O$. }

\bigskip
\noindent {\bf Proof}. Without loss of generality, we may assume that $\psi \geq 0$.
Define
$$ X^0_t\df \cases{X_t
&if $t<\tau $\cr
\partial &if $t\geq \tau,$\cr}
$$
which is reflecting Brownian motion in $D$ killed upon leaving $O$.
It is well-known that $X^0$ is a symmetric Markov process on $\overline D \cap O$
with Dirichlet form $(\EE, W^{1, 2}(D; O^c))$, where
$$  W^{1, 2}(D; O^c)\df \{ u\in W^{1,2}(D): \, u =0 \hbox{ q.e. on } O^c \}.
$$
The transition semigroup for $X^0$ will be denoted by $\{P^0_t, t\geq 0\}$.

Let $O_1$ be a relatively compact open subset of $O$ and let
$f\geq 0$ be $C^1_c$ with $\hbox{supp}[f]\subset O$ and  $f=1$ on
$O_1$. Define $u(x)\df f(x)h(x)$. Then for $x\in D\cap O$,
$$ u(x)-P^0_tu(x)= \E_x \left[ (f(X_0)-f(X_t)) h(X_t); \, t<\tau \right]
+ \E_x\left[ f(X_0) h(X_\tau); \, t\geq \tau \right].
$$
Note that by time-reversal,
$$\eqalignno{
 \int_{D\cap O}& u(x) \E_x \left[ (f(X_0)-f(X_t)) h(X_t); \, t<\tau \right] dx \cr
=& \int_{D\cap O} \E_x \left[ f(X_0)h(X_0)(f(X_0)-f(X_t)) h(X_t); \, t<\tau \right] dx\cr
=& \int_{D\cap O} \E_x \left[ f(X_t)h(X_t)(f(X_t)-f(X_0)) h(X_0); \, t<\tau \right] dx.\cr}
$$
Hence
$$ \eqalignno{
 \int_{D\cap O}& u(x) \E_x \left[ (f(X_0)-f(X_t)) h(X_t); \, t<\tau \right] dx\cr
= &{1\over 2} \int_{D\cap O} \E_x \left[ (f(X_0)-f(X_t))^2 h(X_0) h(X_t); \, t<\tau \right] dx\cr
\leq & {\| h \|^2_\infty \over 2} \int_{D\cap O} \E_x \left[ (f(X_t)-f(X_0))^2 ; \, t<\tau \right] dx.\cr}
$$
Thus
$$\eqalign{
 \limsup_{t\to 0}& {1\over t} \int_{D\cap O} u(x)( u(x)-P^0_tu(x)) dx \cr
 & \leq \limsup_{t\to 0} \Big( {\| h\|^2_\infty \over 2t} \int_{D\cap O}
       \E_x \left[ (f(X_t)-f(X_0))^2 ; \, t<\tau \right] dx\cr
& \qquad +  {\|h\|^2_\infty \over t} \int_{D\cap O}
         f(x)^2 \P_x (t\geq \tau ) dx \Big) \cr
& \leq  \| h \|_\infty^2 \int_{D\cap O} | \nabla f (x) |^2 dx <\infty, \cr}
$$
by Lemma 4.5.2(i) and (4.5.7) of [FOT]. This implies, by Lemma 1.3.4 of [FOT],
 that $u\in W^{1,2}(D; O^c)$ and so $h\in W^{1, 2}(O_1)$.
$$ u(x) = \E_x \left[ u(X^0_{\tau_{O_1}} ) \right] \qquad \hbox{for } x\in O\cap \overline D,
$$
where $\tau_{O_1}\df \inf\{t>0: X^0_t\notin O_1\cap \overline
D\}$. Hence by Theorem 4.3.2 of [FOT], $u$ is $\EE$-orthogonal to
$W^{1,2}(D; O_1^c)$; that is,
$$ {1 \over 2} \int_{D\cap O} \nabla u (x) \cdot \nabla \phi (x)\, dx = 0
\qquad \hbox{for every } \phi \in W^{1,2}(D; O_1^c).
$$
This shows that $u$ and therefore $h$ is harmonic in $D\cap O_1$ with zero Neumann boundary conditions on
$\partial D \cap O_1$. Since $O_1$ is an arbitrary relatively compact open subset of $O$, we conclude
that $u$ is harmonic in $D\cap O$ with zero Neumann boundary conditions on
$\partial D \cap O$.
\qed

\bigskip

The following version of the Neumann boundary Harnack principle is
similar to (but slightly more general) than Theorem 3.9 of [BH].
The result in [BH] is limited to smooth domains whose boundaries
are locally graphs of Lipschitz functions (although the constant
in that theorem depends only on the Lipschitz constant $\lambda$)
and to harmonic functions $h$ as in our Lemma 2.7, with
non-negative $\psi$.

\bigskip

\noindent{\bf Lemma 2.8 (Neumann boundary Harnack principle)}.
{\sl Suppose that $\Phi:\R^{d-1}\to \R$ is a Lipschitz function
with constant $\lambda< \infty$, i.e., $|\Phi( x) -\Phi( y)| \leq
\lambda | x -  y|$ for all $ x,  y \in \R^{d-1}$. Assume that
$\Phi(0)=0$ and let $D=\{x=( x_1, x_2, \dots, x_d)\in \R^d:  x_d>
\Phi(( x_1,\dots,x_{d-1}))\}$. If $r>0$, $c_1>1$, and $h:
B(0,c_1r)\cap D \to [0,\infty)$ is harmonic
with zero Neumann boundary
conditions on $ B(0,c_1r) \cap \prt D$ then
$$h(x) \geq c_2 h(y) \qquad \hbox{for all } x,y \in B(0,r)\cap D,
 \eqno(2.6)
$$
where $c_2>0$ depends only on $\lambda$ and $c_1$.

}

\bigskip
\noindent {\bf Proof}. For $(y_1, \cdots , y_d ) \in \R^d$, denote
$\wt y\df (y_1, \cdots , y_{d-1})$. Define a one-to-one map $\phi:
\phi(\wt y, y_d)=(\wt y, y_d-\Phi(\wt y))$. As $\Phi$ is
Lipschitz, the Jacobians of $\phi$ and its inverse $\phi^{-1}$ are
bounded, with the bound depending only on the Lipschitz constant
$\lambda$. Under $\phi$, ${1\over 2}\Delta $ is mapped into a
uniformly elliptic divergence form operator $L$ with coefficient
matrix $A(x)$ (see Remark 2.1.4 of [K]). Let $U\df \phi(B(0, c_1
r)\cap \overline D)$ and $u(x)\df h(\phi^{-1} (x))$ for $x\in U$.
Using the change of variable formula, we conclude from (2.5) that
for every continuous $\psi \in W^{1,2}(U)$ that vanishes on
$\partial U \cap \{y\in R^d: \, y_d >0\}$,
$$ \int_U A(x) \nabla u (x) \cdot \nabla \psi (x) dx =0. \eqno(2.7)
$$
Let $U^-$ be the ``mirror'' reflection of $U$ with respect to the
hyperplane $\{(\wt y, y_d): y_d=0\}$, that is, $U^-=\{y=(\wt y,
y_d): \, (\wt y, -y_d) \in U\}$. For $y=(\wt y, y_d)\in U_-$,
define $A(y)=A((\wt y, -y_d))$ and $u(y)=u((\wt y, -y_d))$. Then
$L\df \nabla (A \nabla)$ is the uniformly elliptic divergence form
operator defined on the domain $U\cup U^-$. It now follows from
(2.7) and its corresponding version for $U^-$ that
$$  \int_{U\cup U^-} A(x) \nabla u (x) \cdot \nabla \psi (x) dx =0
 \qquad \hbox{for every } \psi \in C_c^\infty (U\cup U^-).
$$
Hence $u$ is a non-negative $L$-harmonic function on $U\cup U^-$.
The desired Harnack inequality for $h$ now follows from the
Harnack inequality for the $L$-harmonic function $u$. \qed

\bigskip

\noindent{\bf Remarks 2.9}. (i) Some regularity conditions for a
harmonic function with zero Neumann boundary conditions have to be
assumed (such as those formulated in Definition 2.6) in order for
the Neumann boundary Harnack principle to hold, even if $D$ has a
$C^\infty$ boundary. The Neumann boundary Harnack principle
does not need to
 hold for a harmonic function in $D\cap B(x_0, r)$ which
satisfies zero Neumann boundary conditions only {\it almost
everywhere} on $\partial D \cap B(x_0, r)$. For example, let $D$ be
a half-space in $\R^d$, $d\geq3$, with $\prt D$ passing through the
origin, and let $h(x) = |x|^{2-d}$. Then $h$ satisfies the Neumann
boundary conditions everywhere except at the origin. The Neumann
boundary Harnack principle does not hold for this function $h$ in
$D\cap B(0,1)$.

(ii) We will apply Lemma 2.8 to two classes of functions. One of
these families consists of harmonic functions defined in a
probabilistic way, as in Lemma 2.7. That lemma shows that Lemma
2.8 is applicable to harmonic functions in this family.

We will also apply Lemma 2.8 to the Green function $x\to G(x,y)$,
where $y\in D\setminus B_*$, and $G(\,\cdot\, ,y)$ is the density
of the expected occupation measure for the reflecting Brownian
motion in a Lipschitz domain $D$ killed upon hitting $B_*$,
starting from $y$. To see that Lemma 2.8 can be applied, consider
any $y \in D\setminus B_*$ and let $U$ be any relatively compact
subdomain of $D\setminus (B_*\cup\{ y\})$. Then for $x\in U$, $
G(x, y)= \E_x [ G( X_{\tau_U}, y)]$. So by Lemma 2.7, $x\to G(x,
y)$ is ``locally'' in $W^{1,2}(D)$ and is harmonic with zero
Neumann boundary conditions on $\partial D$.

\bigskip

\noindent {\bf 3. Simply connected planar domains}.

This section will present some results based on ideas developed in
[BCM], a paper on ``trap'' domains. We will present a new result
on trap domains in Section 5. In this section, we will review only
as much  of the material from [BCM] as is relevant to Problem
1.2. We will use complex analytic notation and concepts. Consult
[Po] for the definitions of prime ends, harmonic measure, etc.

We start with some definitions that apply to domains in any number
of dimensions. Let $X$ be normally reflecting Brownian motion on
$\overline D\subset \R^d$, $d\geq 2$, starting from $x\in D$ and
killed upon hitting a closed ball $B_*$. As is mentioned in the
previous section, $X$ is obtained as the projection of
reflecting Brownian motion $X^*$ on the Martin-Kuramochi
compactification of $D$ into $\overline D$.
The distributions of both $X$ and $X^*$ will be denoted $\P_x$ and
the corresponding expectations will be denoted $\E_x$.
 Let $G(x,y)$ be defined on $(D\setminus B_*)\times (D\setminus B_*)$ by
$$\int_{(D\setminus B_*) \cap A} G(x,y) dy
= \E_x \int_0^{T_{B_*}} {\bf 1}_{\{X_t \in A\}} dt,
 \qquad A\subset \overline D,
$$
\noindent where $dy$ denotes $d$-dimensional Lebesgue measure.
Clearly $G(x, y)$ is a symmetric function on
$(D\setminus B_*)\times (D\setminus B_*)$. It follows from Lemma
3.2 of [CFW] that the function $G(x,y)$ can be extended
continuously to $(D^*\setminus B_*) \times
 (D\setminus B_*)$, where $D^*$ is the Martin-Kuramochi
compactification of $D$ as mentioned in the proof of Theorem 2.2
in the previous section. Note that (cf. Section 2.1 of [BCM])
if $x,y\in D\setminus B_*$ and $x\ne y$ then $t\mapsto G(X^*_t,
y)$ is a continuous local martingale.
 It is easy to see that $G(x,y)$ is the Green function for the
domain $D\setminus B_*$ with (zero) Neumann boundary conditions on
$\prt D$ (in the distributional sense) and (zero) Dirichlet
boundary conditions on $\prt B_*$.

For the rest of this section, suppose that $D$ is a simply
connected open subset of the complex plane $\C$, $z_*$ is the
center of $B_*$, and $z_0$ is a prime end in $D$. Consider a
collection $\{\gamma_n\}_{n\geq 1}$ of non-intersecting cross cuts
of $D$ that do not intersect $B_*$ and such that $\gamma_{n+1}$
separates $\gamma_{n}$ from $z_0$ and the $\gamma_n$'s tend to $z_0$.
Suppose further that $\sigma$ is a curve in $D$ connecting $z_*$
to $z_0$ such that $\sigma \cap \gamma_n$ is a single point $z_n$,
for each $n$. This system of curves divides $D$ into subregions:
let $\Om_n$ denote the component of $D\setminus \gamma_n$ which
does not contain $z_*$. Thus $D_n=\Om_n\setminus\Om_{n+1}$ is the
region between $\gamma_n$ and $\gamma_{n+1}$.  Write
$\Om_1\setminus\sigma = \Om^+ \cup \Om^-$, where each set $\Om^+ $
and $ \Om^-$ is connected, and set $D_n^+=\Om^+\cap D_n$ and
$D_n^-=\Om^-\cap D_n$. See Figure 3.1.

\bigskip
\vbox{ \epsfxsize=3.0in
  \centerline{\epsffile{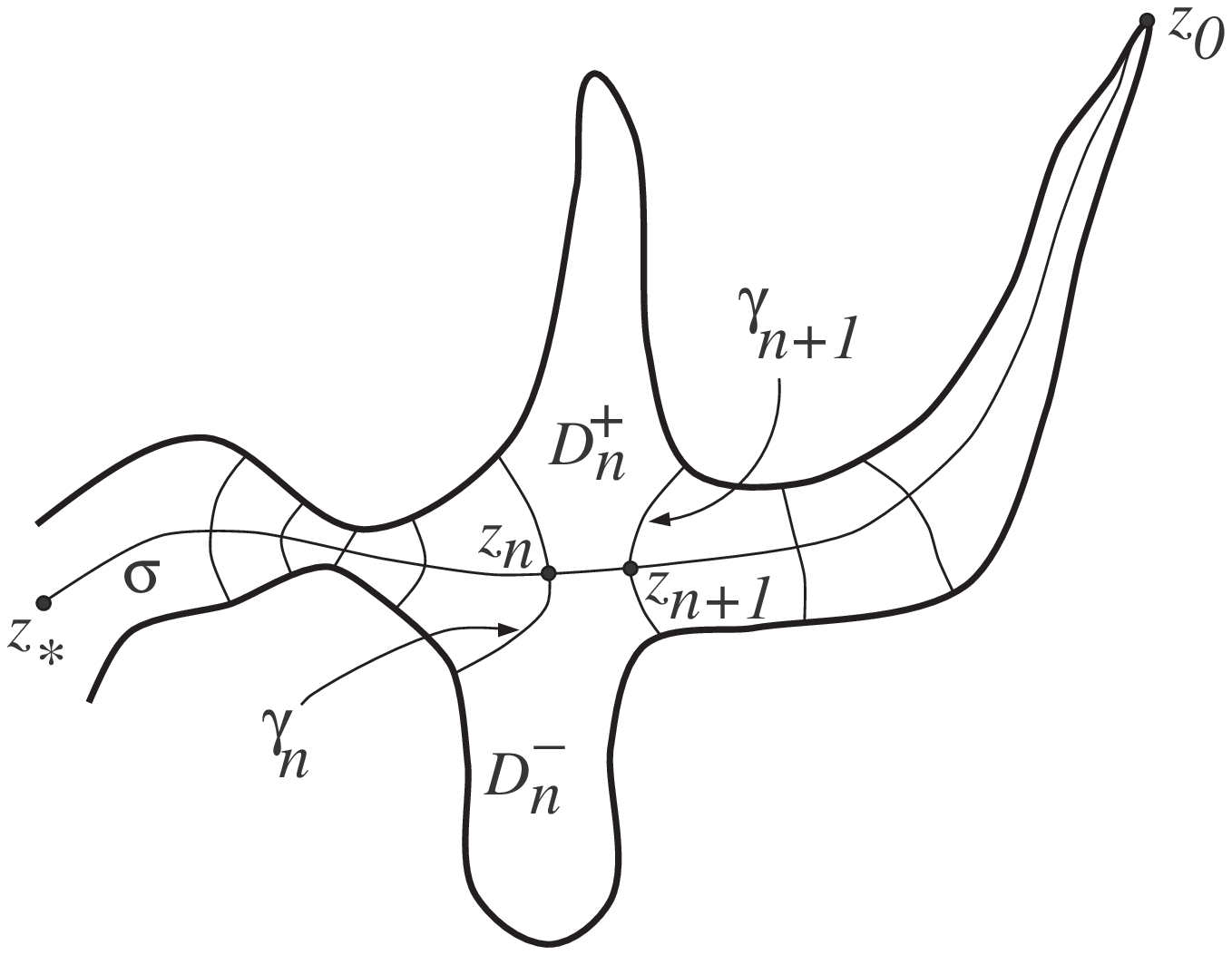}}

\centerline{ Figure 3.1. Hyperbolic blocks.}}
\bigskip

Recall that the harmonic measure of a set $A\subset \b D$ in the
domain $D$, relative to $z$, is denoted $\om(z,A,D)$.

\bigskip
\noindent{\bf Definition 3.1} We will say that the system of
curves $\{\gamma_n\}\cup \sigma$ divide $D$ into {\it hyperbolic
blocks} tending to the prime end $z_0$ if for some $c_* > 0$ and
all $n\geq 1$, the following conditions hold:

\medskip
\item{(i)} $c_* \leq \om(z_*,\b \Om^+ \cap \b D,D)  \leq 1/2$ and
$c_* \leq \om(z_*,\b \Om^- \cap \b D,D)  \leq 1/2$, \item{(ii)}
for all $n\ge 1$ and for all $z\in \b D_n^+\cup\{z_{n-1}\}$, we
have $\om(z,\b D_n^+\cap \b D,D)\ge c_*$, \item{(iii)} for all
$n\ge 1$ and for all $z\in \b D_n^-\cup\{z_{n-1}\}$, we have
$\om(z,\b D_n^-\cap \b D,D)\ge c_*$.

\smallskip

We will call a system of hyperbolic blocks {\it regular} if it
satisfies in addition the following condition,

\item{(iv)} for every $n \geq 1$ there exists $z \in \prt D_n^+
\cap \prt D_n^-$ such that $\om(z,\b D_n^-\cap \b D,D_n)\ge c_*$
and $\om(z,\b D_n^+\cap \b D,D_n)\ge c_*$.

\bigskip

For every simply connected domain and any prime end $z_0$, there
exists a family of regular hyperbolic blocks.
 Here is one way to construct $\{\gamma_n\}_{n\geq 1}$ and $\sigma$.
Suppose that $\varphi$ is a conformal map of the upper half plane
$\H$ onto $D$, such that $\varphi(0) = z_0$ and $\varphi(i) =z_*$.
Then we can take $\gamma_ n = \varphi(\H\cap\{|z|=2^{-n}\})$,
$n\geq 1$, and $\sigma = \{\varphi(iy): 0 < y \le 1\}$. The
conformal invariance of harmonic measure makes it is easy to
verify that the ensemble $\{\gamma_n\}\cup \sigma$ divides $D$ into hyperbolic
blocks tending to $z_0$. Condition (iv) is satisfied by $z =
\varphi( i(3/4)2^{-n})$. Hyperbolic blocks are useful because they
can be constructed geometrically, without knowledge of any
properties of the mapping $\varphi$; see [BCM] for examples of
hyperbolic blocks.

In typical examples, verifying conditions (i)-(iv) is not harder
than verifying just (i)-(iii). We did not include (iv) in the
definition of hyperbolic blocks in order to keep the same
nomenclature as that in [BCM].

\bigskip

\noindent{\bf Theorem 3.2}. {\sl Let $D\in \cD$ be a simply
connected planar domain.

\item{(i)}
If there exist constants $c_*, c'\in (0, \infty)$
such that for each prime end $z_0\in \b D$
there is a system of curves $\{\gamma_n\}\cup
\sigma$ dividing $D$ into
regular
 hyperbolic blocks with parameter $c_*$ and
 $$\sup_{z_0} \sum_n n |\prt D_n\cap \prt D|\le c',\eqno(3.1)$$
then the whole surface of $D$ is active.

\item{(ii)} Let $r_n$ denote the distance between $\gamma_n$ and
$\gamma_{n+1}$.
 If for some prime end $z_0\in \b D$, there is a system of curves
$\{\gamma_n\}\cup \sigma$ dividing $D$ into
regular
 hyperbolic blocks with
$$\sum_{n=1}^\infty n r_n= \infty,    \eqno(3.2) $$
then part of the surface of $D$ is nearly inactive.

\item{} }

\bigskip

Example 3.4 and especially Example 3.6 show that the gap between
parts (i) and (ii) of Theorem 3.2 is not large.

We need a lemma to prove Theorem 3.2.

\bigskip
\noindent{\bf Lemma 3.3} {\sl Suppose that $D\in \cD$ and the
$\{\gamma_k\}$ divide $D$ into regular hyperbolic blocks. For
$n\geq 1$ and $y\in \gamma_{n-1}$, let $x\mapsto h_y(x)$ be the
Poisson kernel with pole at $y$ for reflecting Brownian motion in
$\Omega_{n-1}$ killed upon hitting $\gamma_{n-1}$. Let $\P_x^y$
denote the distribution of Doob's $h_y$-transform of reflecting
Brownian motion on $D^*$ killed upon hitting $\gamma_{n-1}$,
starting from  $x\in \Omega_{n-1}$. Note that $\P_x^y$-a.s., the
process will stay in the closure of $\Omega_{n-1}$ in $D^*$ until
its lifetime. Let $r_n$ denote the distance between $\gamma_{n-1}$
and $\gamma_n$. There exist $c_1,p_1>0$, depending only on $D$ and
$c_*$, such that $\P_x^y (L_{T_{\gamma_{n-1}}}> c_1 r_{n})>p_1$
for any $x\in \gamma_n$.

}

\bigskip
\noindent{\bf Proof}. Let $\varphi$ be a one-to-one conformal map
of $D_{n-1}$ onto the unit disc $S= \{z\in\C: |z| < 1\}$, such
that $\varphi(\prt D_{n-1}^+ \cap \prt D) = I_1\df
\{z=e^{i\theta}: \theta_1 \leq \theta \leq \pi -\theta_1\}$ and
$\varphi(\prt D_{n-1}^- \cap \prt D) = I_2 \df \{z=e^{-i\theta}:
\theta_1 \leq \theta \leq \pi -\theta_1\}$, for some $0 < \theta_1
<\pi/2$. By condition (iv) in Definition 3.1 and conformal
invariance, there exists a point $x_1 \in S$ such that
$\omega(x_1,I_1,S) \geq c_*$ and $\omega(x_1,I_2,S) \geq c_*$.
This easily implies that there exists $\theta_2 =
\theta_2(c_*)\in(\pi/4,\pi/2)$ such that $\theta_1 < \theta_2$.
Let $\theta_3 = (\pi/2 + \theta_2)/2$, $J_1 = \{z=e^{i\theta}:
\theta_3 \leq \theta \leq \pi -\theta_3\}$ and $J_2 =
\{z=e^{-i\theta}: \theta_3 \leq \theta \leq \pi -\theta_3\}$. For
some $c_2=c_2(c_*)>0$ and every $z\in S$ on the imaginary axis,
$\omega(z,J_1\cup J_2,S) \geq c_2$. Let $J^r$ and $J^\ell$ be the
right and left connected components of $\prt S \setminus (I_1 \cup
I_2)$. Let $X$ be  Brownian motion in $S$ with normal reflection
on $I_1 \cup I_2$ killed upon hitting $J^r\cup J^\ell$. It is easy
to see that for some $c_3=c_3(c_*)>0$ and every $z\in J_1\cup
J_2$, if reflecting Brownian motion in $S$ starts from $z$, then
it hits $J^\ell$ before hitting $J^r$ with probability greater
than $c_3$ but less than $1-c_3$.

Let $\gamma_{n-1/2} = \varphi^{-1}(\{z=a+bi\in S: a=0\})$, $K_1 =
\varphi^{-1} (J_1)$ and $K_2 = \varphi^{-1} (J_2)$. By conformal
invariance, for every $x \in \gamma_{n-1/2}$, we have $\omega(z,
K_1\cup K_2, D_{n-1}) \geq c_2$, and for every point $x \in
K_1\cup K_2$, the probability that reflecting Brownian motion
in $D$ starting from $x$ hits
$\gamma_{n-1}$ before hitting $\gamma_n$ is in
the range $(c_3,1-c_3)$.

Find $c_4\in(0,1/8)$ so small that a Brownian motion $W$ starting
from $x$ will make a double loop in an annulus $B(x, r) \setminus
B(x, 3r/4)$ for some $r\in (c_4 r_n, r_n/8)$, and then will make a
crossing from $B(x+ir/2,r/16)$ to the ball $B(x+i2r,r/16)$ within
the convex hull of the two balls, all before leaving $B(x,
r_n/3)$, with probability greater than $1-c_3/2$. A ``double
loop'' in $B(x, r) \setminus B(x, 3r/4)$ means that there exist
$t_1<t_2$ such that $W_t \in B(x, r) \setminus B(x, 3r/4)$ for all
$t\in(t_1, t_2)$, and a continuous version of $t\to \arg(W_t -x)$
increases by $4\pi$ over the interval $[t_1,t_2]$. Note that $c_4$
may be chosen independently of $r_n$, by Brownian scaling.

Consider a reflecting Brownian motion $X_t = x+W_t + N_t$ on $D$,
starting from a point $x \in K_1\cup K_2$. Suppose that $x+W_t$
makes a double loop in an annulus $B(x, r) \setminus B(x, 3r/4)$
for some $r\in (c_4 r_n, r_n/8)$, and then it makes a crossing
from $B(x+ir/2,r/16)$ to the ball $B(x+i2r,r/16)$ within the
convex hull of the two balls, before leaving $B(x,\, r_n/3)$,
during a time interval $[t_1, t_2]$. Suppose moreover, that
$L_{t_2} - L_{t_1} \leq  r /16$. We will show that the two
assumptions taken together yield a contradiction. The second
assumption implies that $|N_t - N_{t_1} | \leq r /16$ for all $t
\in [t_1, t_2]$. This implies that $X$ will make more than one
loop in $B(x, 17r/16) \setminus B(x, 11r/16)$ and then it will
make a crossing from $B(x+ir/2,r/8)$ to the ball $B(x+i2r,r/8)$
within the convex hull of the two balls, before leaving $B(x,\,
r_n/2)$. This is impossible because then $X$ would make a closed
loop around $x$ within $B(x,r_n/2)$, and hence it would have to
cross the boundary of $D$. We conclude that if the first
assumption holds, then $L_{t_2} - L_{t_1} \geq  r /16 \geq c_4 r_n
/16= c_5 r_n$. Since the first event has probability greater than
$1-c_3/2$ and the process $X$ starting from $x \in K_1\cup K_2$
can hit $\gamma_{n-1}$ before $\gamma_n$ with probability greater
than $c_3$, the event that $X$ hits $\gamma_{n-1}$ before
$\gamma_n$ and $L_{t_2} - L_{t_1} \geq c_5 r_n$ has probability
greater than $c_3/2$. This implies that reflecting Brownian motion
in $D$ starting from $x\in \gamma_{n-1/2}$ will hit $\gamma_{n-1}$
before $\gamma_n$ and $L_{T_{\gamma_{n-1}}}\geq c_5 r_n$ with
probability greater than $c_6>0$. Hence, reflecting Brownian
motion in $D$ conditioned to hit $\gamma_{n-1}$ before $\gamma_n$
and starting from $x\in \gamma_{n-1/2}$ will accumulate more than
$c_5 r_n$ units of local time on $\prt D_{n-1}$ before hitting
$\gamma_{n-1}$ with probability greater than $c_6$.

Let $A$ be the interior of $ \ol {D_{n-2} \cup D_{n-1}}$ and let
$\psi$ be a one-to-one conformal mapping of $A$ onto a rectangle
$R=\{a+ib: a_1 < a < a_2, 0 < b < 1\}$, such that $\gamma_{n-2}$
is mapped onto the left side of $R$ and $\gamma_n$ is mapped onto
the right side of $R$. Lemma 3.4 of [BCM] and a simple argument
show that $a_2-a_1$ is bounded above by a constant. Since the
hyperbolic blocks are regular, there exists a point $x\in
\psi(D_{n-2})$ such that the harmonic measure of the upper part of
$R$ in $\psi(D_{n-2})$ is greater than $c_*$, and the same is true
for the lower part of the boundary. An analogous statement is true
for $\psi (D_{n-1})$. All this easily implies that the distance of
$\psi(\gamma_{n-1})$ from the left and right sides of $R$ is
bounded below by $c_7=c_7(c_*)>0$. Let $R_1=\{a+ib: a_1 +c_7/2< a
< a_2-c_7/2, 0 < b < 1\}$. By the Neumann boundary Harnack
principle (Lemma 2.8), for any positive harmonic function $h$ in
$R_1$ with Neumann boundary conditions on the upper and lower
sides of $R_1$, $h(x) \leq c_8 h(z)$ for all $x,z \in
\psi(\gamma_{n-1})$. This applies, in particular, to $h_y\circ
\psi^{-1}$. By conformal invariance, $h_y(x) \leq c_8 h_y(z)$ for
all $x,z \in \gamma_{n-1}$.

Let $g(x)$ be the harmonic function in $D_{n-1}$ with Neumann
boundary conditions on $\prt D_{n-1} \cap \prt D$, equal to $1$ on
$\gamma_{n-1}$ and equal to $0$ on $\gamma_n$. Reflecting
Brownian motion in $D_{n-1}$ conditioned to hit $\gamma_{n-1}$
before $\gamma_n$ is a $g$-transform of the unconditioned process.
We have already proved that the $g$-process starting from $x\in
\gamma_{n-1/2}$ will accumulate more than $c_5 r_n$ units of local
time on $\prt D_{n-1}$ before hitting $\gamma_{n-1}$ with
 probability greater than $c_6$. By the strong Markov property
applied at the hitting time of $\gamma_{n-1/2}$, the same holds if
the starting point belongs to $\gamma_n$.
Without loss of generality, we may and do assume that $h_y(x_0)=1$ for
some $x_0\in \gamma_{n-1}$.
Since $0<c_9<g(x)/h_y(x)<c_{10}<\infty$ for $x\in\gamma_{n-1}$, an
elementary argument shows that the $h_y$-process starting from
$x\in \gamma_{n}$ will accumulate more than $c_5 r_n$ units of
local time on $\prt D_{n-1}$ before hitting $\gamma_{n-1}$ with
probability greater than $c_6$.
 \qed

\bigskip
\noindent{\bf Proof of Theorem 3.2}. (i) Let $d_D(x)\df \dist(x,
\prt D)$. Consider $x_0 \in D$. It is not hard to see that there
exists $z_0\in \prt D$ and a corresponding family of $\gamma_n$'s
such that $x_0 \in D_{n_0}$ for some $n_0$ and $\dist(x_0, \prt
D_n) \geq c_1 d_D(x_0)$, where $c_1\in (0, 1)$ is a constant
depending only on $D$. By the proof of Theorem 2.2 (see especially
Lemmas 3.4 and 3.5) in [BCM], $G(x_0, \,\cdot \,)$ is bounded by
$c_2 k$ on $D_{k}$ for $k\leq n_0-1$. Hence $G(x_0, \,\cdot \,)$
is bounded by $c_2 n_0$ on $D_{n_0-1}$. By the Harnack principle,
it is bounded by $c_3n_0$ on $\prt B(x_0, c_1 d_D(x_0)/2)$, and
since
$$
G(x_0, x)=\E_x \left[ G(x_0, X_{T_{B(x_0,  c_1 d_D(x_0)/2)}}) \right]
\qquad \hbox{for }  x\in D\setminus B(x_0,  c_1 d_D(x_0)/2),
$$
the same bound holds on $D \setminus B(x_0, c_1 d_D(x_0)/2)$.
We obtain,
 $$\eqalign{
 \E_{x_0} L_{T_{B_*}}
 &= \sum_{n=1}^\infty \E_{x_0}\int_0^{T_{B_*}}
 \bone_{\prt D_n \cap \prt D}(X_t)
 dL_t\cr
 &= \sum_{n=1}^{n_0 -1}
 \E_{x_0}\int_0^{T_{B_*}} \bone_{\prt D_n \cap \prt D}(X_t) dL_t
 +\sum_{n=n_0}^\infty \E_{x_0}\int_0^{T_{B_*}}
 \bone_{\prt D_n \cap \prt D}(X_t)  dL_t\cr
 &= \sum_{n=1}^{n_0 -1}
 \int_{\prt D_n \cap \prt D} G(x_0, x) \sigma (dx)
 +\sum_{n= n_0}^\infty  \int_{\prt D_n \cap \prt D} G(x_0, x) \sigma (dx)\cr
 &\leq  \sum_{n=1}^{ n_0-1} |\prt D_n \cap \prt D|
 \sup_{x\in \prt D_n \cap \prt D} G(x_0,x)
 + \sum_{n=n_0}^\infty |\prt D_n \cap \prt D|
 \sup_{x\in \prt D_n \cap \prt D} G(x_0,x)\cr
 &\leq  \sum_{n=1}^{ n_0-1} |\prt D_n \cap \prt D|\,  c_2 n
 + \sum_{n=n_0}^\infty |\prt D_n \cap \prt D|\,  c_3 n_0\cr
 &\leq  \sum_{n=1}^\infty c_4 \, n\, |\prt D_n \cap \prt D| .
}$$
This is bounded by a constant independent of $x_0$, by assumption
(3.1). Hence we obtain $\sup_{x\in D} \E_{x} L_{T_{B_*}} < \infty$
and, therefore, $ \inf_{x\in D} u(x) =\inf_{x\in D}
\E_x \exp(-L_{T_{B_*}}) >0$. This means that the whole surface of
$D$ is active.

(ii) Find a prime end $z_0 \in \prt D$ and a family of
$\gamma_n$'s such that (3.2) holds, that is,
$$\sum_{n=1}^\infty n r_n=\infty.
$$
Note that in the case of a simply connected domain $D$ in $\R^2$,
the Martin-Kuramochi boundary $D^*$ of $D$ and the corresponding
reflecting Brownian motion $X^*$ on $D^*$ can be realized as
follows. Let $\varphi$ be a conformal map from $\H\df \{a+bi: \,
b>0\}$ to $D$ and define $D^*$ to be the union of $D$ and its
prime ends. The map $\varphi$ extends continuously to a one-to-one map from
$\overline \H$ to $D^*$. Let $Y$ be reflecting Brownian motion
on $\overline \H$. Then $\varphi (Y)$ is a time change of
reflecting Brownian motion $X^*$ on $D^*$. We will use this
constructed reflecting Brownian motion $X^*$ in this proof. Recall
that $G(z_0, x)$ is well defined for $x\in D^*\setminus \{z_0\}$
by the second paragraph of this section. For $ a\geq 0$, define
$$ \eta_a \df  \left\{x\in D^* \setminus (B_*\cup \{z_0\}): G(z_0, x) = a
   \right\}.
$$
First, we claim that there exist positive integers $m_0$, $m_1$
and a positive constant $a_0$ such that there is at least one
$D_n$, but at most $m_1$
 such sets, between $\eta_a$ and $\eta_{a+m_0}$, for every $a>a_0$.

Recall that $z_*$ is the center of $B_*$. Let $\varphi: \H \to D$
be a one-to-one conformal mapping, such that $\varphi(0)=z_0$ and
$\varphi(i) =z_*$. Define $h:\, D \to \R$ by $h(z) = -\log
|\varphi^{-1}(z)|$. Then $h$ is harmonic in $D$ with Neumann
boundary conditions and a pole at $z_0$. Let $\eta^*_a= \{x \in D:
h(x) = a\}$. Lemma 3.4 of [BCM] and conformal invariance easily
imply that there exists an integer $\widetilde m_0$ such that for
any $a\in \R$ there is at least one $D_n$ between $\eta^*_a$ and
$\eta^*_{a+ \widetilde m_0}$. It follows from the conformal
invariance of the Green function that $h_1(z)\df G(z_0, \varphi
(z))$ is the Green function for reflecting Brownian motion in $\H$
starting from $0$ and killed upon hitting $\varphi (B_*)$. It is
easy to see that $h_1(z)$ and $-\log |z|$ are comparable on $\H
\cap \{z: |z|<r\}$, for some $r>0$. This implies the existence of
a positive integer $m_0$ and a constant $a_0>0$ such that there is
at least one $D_{n}$ between $\eta_a$ and $\eta_{a+m_0}$ for every
$a>a_0$. From (3.1), the inequalities preceding (3.2) and (3.3) in
[BCM] as well as Lemma 3.5 of [BCM], we see that there exists
$m_1<\infty$ such that there are at most $m_1$ sets $D_n$ between
any $\eta_a$ and $\eta_{a+m_0}$.

Let $\alpha_j$ be the sum of $nr_n$ restricted to
integers $n$ such that $D_{n-1}$ lies between $\eta_a$ and
$\eta_{a+m_0+1}$, where $a$ is of the form $km_0 +j$. Every set
$D_{n-1}$ lies between $\eta_a$ and $\eta_{a+m_0+1}$ for some
integer $a$, namely for the largest integer $a$ such that
$\eta_a\cap \Omega_{n-1} =\emptyset$. This and (3.2) imply that
$\sum_{j=0}^{m_0} \alpha_j =\infty$. We will assume without loss
of generality that $\alpha_0=\infty$.

We define $k(n)$ to be the integer $k$ which maximizes $kr_k$
among all $k$'s such that $D_{k-1}$ lies between $\eta_{(n-1)m_0}$
and $\eta_{nm_0}$ (we take the largest of the $k$'s with these
properties if the above definition does not uniquely identify
$k(n)$).
 If we restrict the sum in (3.2) to $k(n)$'s, its value
will be infinite, because there are at most $m_1$ sets $D_{n-1}$
between any $\eta_a$ and $\eta_{a+m_0}$. By Lemma 3.5 of [BCM] and
the comparability of $-\log |\varphi (z)|$ and $G(z_0, z)$ for $z$
in a neighborhood of $z_0$, $c_1 n \leq k(n)\leq c_2 n$.

By (3.1), the inequalities
preceding
 (3.2) and (3.3) in [BCM] as well as Lemma 3.5 of [BCM], $c_1 k
\leq G(z_0,x)\leq c_2 k$ for $x\in D_{k}$ for $k\geq 1 $. Let
$\beta $ be the smallest integer multiple of $m_0$ greater than
$\max \{2, \, (c_2/c_1)\}$. We have
 $$
\sum_{j=1}^\infty\,  \sum_{ n=\beta^{2j-1+m} }^{\beta^{2j+m} } nr_{k(n)}
 =\infty
$$
for $m=0$ or $1$ and we will assume without loss of generality
that we can take $m=0$, i.e.,
 $$\sum_{j=1}^\infty \,
 \sum_{n=\beta^{2j-1}}^{\beta^{2j}} nr_{k(n)}=\infty.\eqno(3.3)$$

Let  $X^*$ be reflecting Brownian motion on $D^*$ starting from
some $x_0\in D_{n_0}$, where $n_0$ is large. Define
 $$\eqalign{
 S_j &\df \inf\{t>0: X^*_t \in \eta_{\beta^{2j}}\}, \quad j\geq 1,\cr
 T^{j,n}_1 &\df  \inf\{t>S_j: X^*_t \in\eta_{nm_0}\}, \quad n\geq 1, \cr
 U^{j,n}_k &\df \inf\{t>T^{j,n}_k: X^*_t \in\eta_{(n-1)m_0}\}, \quad
 n,k\geq 1, \cr
 T^{j,n}_k &\df \inf\{t>U^{j,n}_{k-1}: X^*_t \in\eta_{nm_0}\}, \quad
 n,k\geq 2, \cr
 N^j_n &\df \max \{ k: U^{j,n}_k \leq S_{j-1}\},\quad n\geq 1.
 }$$
In other words, $N^j_n$ is the number of downcrossings of
$[(n-1)m_0,nm_0]$ by $M_t\df G(z_0, X^*_t)$ between times $S_j$
and $S_{j-1}$. This is of interest to us only for $n$ such that
$[(n-1)m_0,nm_0] \subset [\beta^{2j-2},\beta^{2j}]$. The process
$M$ is a continuous local martingale so it is a time-change of
Brownian motion, until it hits $0$.

Consider a one-dimensional Brownian motion $W$ starting from
$\beta^{2j}$ and killed at the hitting time $T$ of $\beta^{2j-2}$.
It follows easily from the Ray-Knight Theorem that there is an
event $A$ with probability greater than $p_1>0$, such that on $A$,
the local time $L^x_T$ accumulated by $W$ at the level $x$ before
time $T$ is greater than $c_4 \beta^{2j}$ for all $x\in
(\beta^{2j-1},\beta^{2j})$. We will apply excursion theory to
excursions of $W$ from the set $\{nm_0: \, \beta^{2j-1} \leq nm_0
\leq \beta^{2j}\}$. Given the local time $\{L^x_T, x=nm_0\in
[\beta^{2j-1},\beta^{2j}]\}$ and assuming the event $A$ occurs,
the distribution of the number of excursions going from $nm_0$ to
$(n-1)m_0$ is minorized by a Poisson random variable with
expectation $K_n \geq c_5 \beta^{2j} / m_0\df c_6 \beta^{2j}$.
Conditional on $\{L^x_T, x=nm_0\in [\beta^{2j-1},\beta^{2j}]\}$,
these random variables are independent. Let $T^M(b) = \inf\{t>0:
M_t = b\}$ and $M^j_t = \{M_{t+ T^M(\beta^{2j})}, t\in [0,
T^M(\beta^{2j-2}) - T^M(\beta^{2j})]\}$. Since $M^j_t$ is a
time-change of $W_t$, there exists an event $A'$ with
$\P_{x_0}(A')>p_1$, such that on $A'$, conditional on the local
time of $M^j$, the numbers of excursions of $M^j$ between
consecutive points of $\{nm_0,\beta^{2j-1} \leq nm_0 \leq
\beta^{2j}\}$ are independent random variables minorized by
independent Poisson random variables with means $K_n\geq c_6
\beta^{2j}$.

Note that the processes $M^j$ are independent. We will now
condition the process $X^*$ on the local times of $M^j$'s and the
endpoints of excursions of $X$ from $\{\eta_{nm_0}, \beta^{2j-1}
\leq nm_0 \leq \beta^{2j}\}$.

Recall that $k(n)$ is an integer such that $D_{n-1}$ lies between
$\eta_{(n-1)m_0}$ and $\eta_{nm_0}$. An easy argument based on
Lemma 3.3 shows that given endpoints of an excursion of $X^*$
going from $\eta_{nm_0}$ to $\eta_{(n-1)m_0}$, the amount of local
time accumulated by the excursion on $\prt D$ is greater than $c_7
r_{k(n)}$ with probability greater than $p_2>0$.

Let $J_n$ be the distribution of the local time accumulated by
$X^*$ on the part of $\prt D$ between $\eta_{(n-1)m_0}$ and
$\eta_{nm_0}$, during the time interval
$(T^M(\beta^{2j}),T^M(\beta^{2j-2}))$. We have shown that on an
event $A_j$ of probability greater than $p_1$, $J_n$ is
stochastically minorized by a random variable $I_n$ whose
distribution is Poisson with mean greater than $p_2 c_6 \beta^{2j}
\cdot c_7 r_{k(n)}$. Hence $J_n$ is minorized by a random variable
$I_n$ which has mean $\lambda_n$ greater than $ c_8
r_{k(n)}\beta^{2j}$ and variance $\lambda_n$. Moreover, we can
assume that the $I_n$'s are independent given $A_j$. Hence, the local
time accumulated by $X$ between hitting of $\eta_{\beta^{2j}}$ and
$\eta_{\beta^{2j-2}}$, on the part of $\prt D$ between these
curves, is stochastically minorized by a random variable $H_j$
such that on the event $A_j$, its mean is bounded below by
$\sum_{j: \, \beta^{2j-1} \leq nm_0 \leq \beta^{2j}}c_8
r_{k(n)}\beta^{2j}\geq \sum_{j:\, \beta^{2j-1} \leq nm_0 \leq
\beta^{2j}}c_9 n r_{k(n)}$ and the variance is equal to its mean.
It follows that $H_j$ takes a value larger than $b_j \df {1\over
2} \sum_{j: \, \beta^{2j-1} \leq nm_0 \leq \beta^{2j}}c_9 n
r_{k(n)}$ with probability greater than $p_2>0$. Since the $M^j$'s
are independent, we can assume that the $H_j$'s are independent.
Let $\Lambda_j$ be independent random variables with $P(\Lambda_j
=b_j) = 1- P(\Lambda_j = 0) = p_2$. Since the reflecting Brownian
motion $X^*$ starting from $x_0\in D_{n_0}$ has to go through
$\gamma_j$ for $j=n_0, n_0-1, \cdots , 1$ before reaching
$\gamma_0\df \prt B_*$, the distribution of the local time
accumulated by $X^*$ before hitting $B_*$ is minorized by the
distribution of $\sum_{j=1}^{n_0} \Lambda_j$. In view of (3.3),
$\sum_{j= 1}^\infty b_j = \infty$, and this easily implies that
$\sum_{j=1}^\infty \Lambda_j=\infty$, a.s. Hence, for any $b\in
(0,  \infty)$, there is some $n_0$ such that $\P \left(
\sum_{j=1}^{n_0} \Lambda_j>b\right)> 1-1/b$. This implies that for
any $x_0\in D_{n_0}$, $\P_{x_0}(L_{T_{B_*}}>b)> 1-1/b$. Therefore,
$\inf_{x\in D}\E_x\exp(-L_{T_{B_*}})=0 $ and we see that part of
the  surface of $D$ is nearly inactive.
 \qed

\bigskip

\noindent{\bf Example 3.4}. Our first example is very simple.
Suppose that for some $\alpha>1$,
 $$
D= \left\{x = (x_1, x_2): \, | x_2 | \leq x_1^{ \alpha}
\hbox{ and }  0 < x_1 < 1 \right\}.
 $$
The interesting range of the parameter is $\alpha>1$. We will show
that if $\alpha \in (1,2)$ then the whole surface of $D$ is active
and when $\alpha \geq 2$ then it is not.

It is easy to see that it is sufficient to analyze only one
boundary point, namely, $(0,0)$. We generate a corresponding system
of hyperbolic blocks by letting $\gamma_n$'s be vertical cuts of
the domain at distance $2^{-k} + j 2^{-k\alpha}$ from 0, for
all $j\geq 0$ such that $2^{-k} + j 2^{-k\alpha}\leq 2^{-k+1} -
2^{-k\alpha}$, for all $k\geq 2$.

The number of hyperbolic blocks whose distance from 0 is between
$2^{-k}$ and $2^{-k+1}$ is of order $2^{-k(1-\alpha)}$. Hence the
blocks in this family have indices $n$ of order $\sum_{j\leq k}
2^{-j(1-\alpha)} \approx 2^{-k(1-\alpha)}$. The perimeter of each
of these blocks is of order $2^{-k \alpha}$, so the contribution
from these blocks to the sum in (3.1) is of order
$2^{-k(1-\alpha)}\cdot 2^{-k(1-\alpha)} \cdot 2^{-k \alpha} =
2^{-k (2-\alpha)}$. If $\alpha < 2$ then $\sum_{k\geq 1} 2^{-k
(2-\alpha)}< \infty$, so part (i) of Theorem 3.2 implies that the
whole surface of $D$ is active.

The distance between $\gamma_n$ and $\gamma_{n+1}$ is comparable
to the perimeter of $D_n$, so the same calculation as above shows
that the sum in (3.2) is comparable to $\sum_{k\geq 1} 2^{-k
(2-\alpha)}$ and this is infinite for $\alpha \geq 2$. Therefore
part of $\partial D$ is nearly inactive when $\alpha \geq 2$.

The multidimensional version of this example will be discussed
in Example 4.13.
\bigskip

It is interesting to compare the above result with the
semimartingale property of reflecting Brownian motion $X$ in $D$
starting from the tip ${\bf 0}\df (0, 0)$. It is shown in
DeBlassie and Toby [DT] that $X$ starting from ${\bf 0}$ is a
semimartingale if and only if $\alpha <2$. See also Theorem 3.1(i)
of Burdzy and Toby [BT] for a similar result. Fukushima and
Tomisaki [FT] proved for domains in the shape of multidimensional
cusps that reflecting Brownian motion starting from the cusp point
is a semimartingale if $\alpha <2$. We will show in Remark 4.14
below that it is not a semimartingale when $\alpha \geq 2$. \qed

\bigskip
\noindent{\bf Remark 3.5}. In the definition of $\cD$, it is
required that $| \prt_L D|$ be finite. One can of course relax
this condition  using localization. However if $| \prt_L
D|=\infty$ (under whatever generalization one uses) and if $u$ is
a weak solution to (1.3)-(1.5) in the sense of (2.3) with $0\leq
u\leq c$, then $\inf_{x\in \partial_L D} u(x)=0$. For suppose
otherwise, that is, there exists $c_0>0$ such that $u(x)\geq c_0$
for every $x\in \prt_L D$. Let $g$ be a smooth function with
compact support in $\R^d$ such that $g=0$ on $B_*$ and $g=1$ on
$\partial D$. Then by (2.3) we have
$$ \int_{D\setminus B_*} \nabla g (x)\cdot \nabla u (x)  dx
= -c \int_{\prt_L D} g(x) u(x) \sigma (dx) =-\infty.
$$
This is impossible since the left hand side should be finite by
the Cauchy-Schwarz inequality. \qed

\bigskip
\noindent{\bf Example 3.6}. We will analyze a fractal domain which
contains channels that become thinner at the same rate as the
single channel in Example 3.4. In the present example, the
distance between $\gamma_n$ and $\gamma_{n+1}$ is much smaller
than the perimeter of $D_n$ for some $n$. Nevertheless, there is
no gap between conditions (3.1) and (3.2) for this family of
domains.

Suppose that $\alpha>0$, $\beta >1$ and let $a_k = \sum_{j=1}^k
2^{-(j-1)\alpha}$. Let $\cS_n$ be the family of all binary
(zero-one) sequences of length $n$. We will write $\s = (s_1, s_2,
\dots, s_n)$ for $\s\in \cS_n$. For integer $k\geq 1$ and $\s\in
\cS_k$, we set $b_{\s} = \sum_{j=1}^k s_j 2^{-j}$. Let $A_* =
[0,1]^2$, for $k\geq 1$ and $\s\in \cS_k$ let
 $$A_{\s} = \{(x_1,x_2)\in \R^2: a_k \leq x_1\leq a_{k+1},
 b_{\s}\leq x_2 \leq b_{\s} + 2^{-k\beta}\},$$
and let $D$ be the connected component of the interior of $A_*\cup
\bigcup_{k\geq 1} \bigcup_{\s\in\cS_k} A_{\s}$ that contains the
open square $(0, 1)^2$ (see Figure 3.2). (Note that when $\beta
\geq 2$, the interior of $A_*\cup \bigcup_{k\geq 1}
\bigcup_{\s\in\cS_k} A_{\s}$ is disconnected.)

\bigskip
\vbox{ \epsfxsize=3.5in
 \centerline{\epsffile{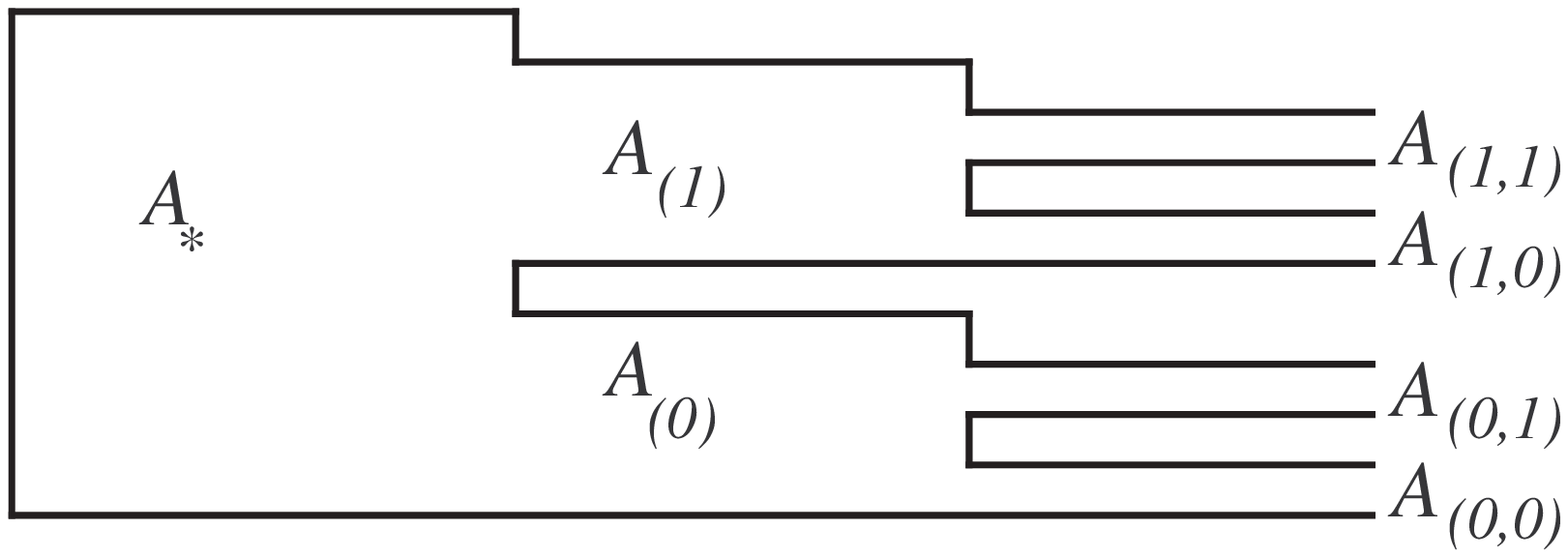}}

\centerline{Figure 3.2.} }
\bigskip

The interesting range of parameters is $\beta >\alpha>0$. If
$\alpha \leq 1$ then $|\prt D|=\infty$, so part of the  surface
of $D$ is nearly inactive by Remark 3.5. We will show that if $1<\alpha <
\beta < 2\alpha$, then the whole surface of $D$ is active and when
$\beta \geq 2\alpha>2$, then part of the surface is nearly inactive.

As in the case of Example 3.4, we will analyze only the family of
hyperbolic blocks corresponding to a boundary point at the end of
a channel. The analysis of other boundary points is
straightforward but tedious so it is omitted. Fix a boundary point
$z_0$ at the end of a channel, i.e., a point whose first
coordinate is $\sum_{j=1}^\infty 2^{-(j-1)\alpha}$. Let ${\cal
A}_k$ be the family of vertical lines $K_{k,n}=\{(x,y): x=a_k + n
2^{-k \beta}\}$, with $n\geq 1$ such that $a_k + n 2^{-k \beta}
\leq a_{k+1} $. Let ${\cal C}_k$ be the family of these line
segments in $K_{k,n} \cap D$ for $K_{k,n} \in {\cal A}_k$ that
separate $z_0$ from $A_*$. Let $\{\gamma_n\}$ be the relabelled
family $\bigcup_k {\cal C}_k$.

Let $\{\s_k\in \cS_k, k\geq 1\}$ be the sequence such that the
channel formed by the $A_{\s_k}$'s approaches $z_0$.
 The number of hyperbolic blocks $D_n$ defined by the $\gamma_n$'s
needed to reach $A_{\s_k}$ is of order $\sum_{j\leq k}
2^{-j\alpha}/2^{-j \beta}\approx 2^{k(\beta-\alpha)}$. Consider a
hyperbolic block $D_n$ which intersects $A_{\s_k}$ The set $D_n$
may be either a square or it may contain a ``tree'' of thin
channels. Consider first $D_n$'s that are squares. There are about
$2^{k(\beta-\alpha)}$ such hyperbolic blocks, so they correspond
to $n$ in (3.1) of order $\sum_{j\leq k} 2^{j(\beta-\alpha)}$,
which is within a constant multiple of $2^{k(\beta-\alpha)}$. The
perimeter of any such $D_n$ is of order $2^{-k\beta}$, so the
total contribution of such $D_n$'s to (3.1) is of order
$2^{k(\beta-\alpha)}2^{k(\beta-\alpha)}2^{-k\beta} \approx
2^{k(\beta -2\alpha)}$. The series $\sum_k 2^{k(\beta -2\alpha)}$
is summable if and only if $\beta < 2 \alpha$.

Next consider a $D_n$ which intersects $A_\s$ with $\s\in\cS_k$
and contains a side ``tree'' of thin channels. The length of its
boundary is of order $ \sum_{j\geq k} 2^{(j-k)} 2^{-j\alpha}
\approx 2^{-k\alpha}$. It corresponds to $n$ in (3.1) of order
$2^{k(\beta-\alpha)}$. There are at most two such $D_n$'s for each
$A_\s$, so their contribution to (3.1) is of order
$2^{k(\beta-\alpha)}2^{-k\alpha} \approx 2^{k(\beta -2\alpha)}$.
Hence, the contribution of $D_n$'s with side channels is of the
same order as the contribution of $D_n$ that have the square
shape. We conclude that (3.1) holds if $\beta < 2 \alpha$.

If $\beta \geq 2 \alpha$, then the contribution of the square
$D_n$'s is enough to make the left hand side of (3.2) infinite,
due to the estimates presented above.
\qed

\bigskip

\noindent {\bf 4. Higher dimensional domains}.

This section is devoted to a family of multidimensional domains.
The family may seem small, but it contains many examples that
arise naturally in the context of the present paper and that of
[BCM].
Before presenting the main result of this section, Theorem 4.3, we
will state a definition and a technical assumption.

Recall the family $\cD$ from Definition 2.1.

\bigskip

\noindent{\bf Definition 4.1}. We will say that a domain $D\subset
\R^d$, $d\geq 3$, belongs to the family $\cD_1$ if $|D|<\infty$,
$D\in \cD$, and it satisfies the properties listed below. Recall
the meaning of the Lipschitz constant $\lambda$ from Definition
2.1. Let $\gamma_0\df \prt B_*$. For every boundary point $z \in
\prt D$ there exists a family of disjoint smooth
$(d-1)$-dimensional surfaces $\{\gamma_n\}_{n\geq 1}$, such that
$\gamma_n\subset D$, and the set $D\setminus \gamma_n$ consists of
two open connected components, $\Omega_n$ and $\Omega_n'$. For
every $n$, we have $z\in \ol\Omega_n$, $B_* \subset \Omega'_n$,
and $\Omega_{n+1} \subset \Omega_n$. Let $r_n$ be the distance
between $\gamma_n$ and $\gamma_{n+1}$. There exist $k_0<\infty$
and $0<\alpha_1, \, \alpha_2,\dots, \alpha_7<\infty$, depending
only on $D$, satisfying the following conditions.

\item{(i)} For $n\geq 0$, $\alpha_1 <r_n/r_{n+1} < \alpha_2$ and
$\gamma_n$ can be covered by at most $k_0$ balls of radius
$\alpha_3 r_n$.

\item{(ii)} For every $n\geq 0$, there exists a curve
$\Gamma\subset D$ of length less than $\alpha_4 r_n $, connecting
$\gamma_n$ and $\gamma_{n+1}$, whose distance from $\prt D$ is
greater than $\alpha_5 r_n$.

\item{(iii)} For every $n\geq 0$ and $x\in \ol\gamma_n\cap\prt D$,
there exists an orthonormal coordinate system $CS$ with the
property that $\prt D \cap  B (x,\alpha_6 r_n)$ is the graph of a
Lipschitz function with constant $\lambda$ in $CS$.

\item{(iv)} For every $n\geq 0$ and $x\in \ol\gamma_n\cap\prt D$,
there exists an orthonormal coordinate system $CS$ with the
property that $\prt \Omega_n \cap  B (x, \alpha_6 r_n)$ is the
graph of a Lipschitz function with Lipschitz constant $\alpha_7 $
in $CS$, and the analogous statement is true for $\Omega'_n$ in
place of $\Omega_n$.

We will write $D_n \df \Omega_n \setminus \overline{\Omega_{n+1}}$.

\bigskip

Note that it follows from part (ii) of Definition 4.1 that there
is a constant $\alpha_{8}>0$ depending only on $D\in \cD_1$ such
that $|D_n|\geq \alpha_{8}r_n^d$ for every $n\geq 0$.

Our proof of the second part of our main result in this section,
Theorem 4.3, requires the following technical assumption, Condition 4.2.
We will discuss ways of verifying this assumption
after the proof of Theorem 4.3.

\bigskip

\noindent {\bf Condition 4.2.} {\sl There exist $0<m_0\leq m_1<\infty$ such that
for any $z\in\prt D$ and $\gamma_n$'s as in Definition 4.1, if
$n>m_1$ and $x_0\in \Omega_n$, then $\sup_{x\in \gamma_{n-m_0}}
G(x_0,x) \leq \inf _{x\in \gamma_n} G(x_0,x)$. Here $G(x, y)$ is
the Green function of reflecting Brownian motion in $D$ killed
upon hitting $B_*$.}

\bigskip
Recall that we say that the whole surface of $D$ is active if
(1.6) holds.

\bigskip

\noindent{\bf Theorem 4.3}. {\sl Suppose that $D\subset \R^d$,
$d\geq 3$,
is
such that $D\in\cD_1$.

(i) If for each boundary point $z\in \b D$, there exists a system
of surfaces $\{\gamma_n\}$ as in Definition 4.1 such that
 $$\sup_{z\in \prt D} \sum_{n\geq 1} |\prt D_n \cap \prt D|
 \sum_{k=1}^n r_k^{2-d}< \infty,\eqno(4.1)$$
then the whole surface  of $D$ is active. Here $|\prt
D_n \cap \prt D|$ denotes the $(d-1)$-dimensional surface measure
of $\prt D_n \cap \prt D$.

(ii) Suppose now that Condition  4.2 holds. If there exists
a boundary point $z\in \prt D$ and a family of surfaces
$\{\gamma_n\}$ as in Definition 4.1, such that
 $$ \sum_{n\geq 1} r_n^{d-1}
 \sum_{k=1}^n r_k^{2-d}=\infty,\eqno(4.2)$$
then part of the surface of $D$ is nearly inactive.

}
\bigskip

The proof of the above theorem will be preceded by a few lemmas.
Recall from Section 1 that $B_*\subset D$ is a fixed reference
ball. In our proofs, $c_j, k_j, m_j$ and $p_j$, $j=0,1,\dots$,
will denote strictly positive and finite constants depending only
on $D$.

\bigskip

\noindent{\bf Lemma 4.4}. {\sl Let $D\in
\cD_1$, $z_0 \in \prt D$, and let $\{\gamma_n\}$ and $\{r_n\}$ be
defined relative to $z_0$ as in Definition 4.1. There exist $c_1,
c_2 \in (0,\infty)$, depending only on $D$, such that
 $$c_1 \sum_{k=1}^n r_k^{2-d}
 \leq G(x, y) \leq c_2 \sum_{k=1}^n r_k^{2-d},$$
for all $n\geq 1$, $x\in \gamma_n$ and $y\in \gamma_{n+1}$.
}

\bigskip

\pf
Let $\{\gamma_n, n\geq 0\}$ and $\{\Omega_n, n\geq 0\}$ be as in
Definition 4.1, and recall that $G(x, y)$ is the
 Green function for reflecting Brownian motion $X^*$ on $D^*$
killed upon hitting $B_*$. As we observed in Section 2, $X^*=X$ on
$D\cup \prt_L D$. For $k\geq 0$, let
$$G_{\Omega_k} (x, y)\df G(x, y)-
\E_x \left[ G(X_{T_{\gamma_k}}, y) \right], \qquad x, y \in \Omega_k,
$$
and note that $G_{\Omega_k} (x, y)$
 is the Green function for reflecting Brownian motion in
$\Omega_k$ killed upon hitting $\gamma_k$. Note
that
 since $\gamma_0 \df \partial B_*$, $G_{\Omega_0}(x, y) =G(x, y)$.

It follows easily from Definition 4.1 that we can find points
$z_k\in D_k$, $k\geq 0$, and finite positive constants $c_0<c_1$
depending only on $D$ such that
 $$ \eqalign{\hbox{dist}\,(z_k, \partial D_k)> c_0 r_k \ &\hbox{ and } \
 \max\left\{\hbox{dist}\,(z_k, \gamma_k),
 \  \hbox{dist}\,(z_k, \gamma_{k+1})\right\} < c_1r_k\cr
 & \hbox{for every } k\geq 0.\cr}
 $$
Let $B_k\df B(z_k, c_0 r_k/2)$. Starting at any point
in
 $\partial B(z_k, c_0 r_k/4)$, the expected time that Brownian
motion
spends in $B_k$ before hitting $\partial D_k$ is larger than $c_2
r_k^2$. By the support theorem for standard $d$-dimensional
Brownian motion, starting from any point
in
$\partial B(z_k, 3c_0 r_k/4)$, there is probability at least
$p_1>0$ (not depending on $k$) that the Brownian motion will hit
the ball $B(z_k, c_0 r_k/4)$ before hitting $\partial D_k$.  So
starting at such a point the expected time spent in $B_k$ before
hitting $\partial D_k$ is at least $p_1 c_2 r_k^2$. This
implies that
$$ \int_{B_k} G_{\Omega_k} (x, y) dx \geq p_1 c_2 r_k^2
 \qquad \hbox{for every }
  y\in \partial B(z_k, 3c_0r_k/4).
$$
Using the Harnack inequality and the fact that $|B_k|=c_3 r_k^d$, it follows that
$$ G_{\Omega_k}(z_k,y)\geq c_4r_k^{2-d} \qquad \hbox{for every } y\in \partial B(z_k, 3c_0r_k/4).
$$
By the Harnack and the Neumann boundary Harnack principle (Lemma 2.8), we have
$$G_{\Omega_k}(x,y)\geq c_5 r_k^{2-d}, \qquad x\in \gamma_{k+1} \hbox{ and }
   y\in \gamma_{k+2}. \eqno (4.3)$$

On the other hand, starting in $B_k$ the expected amount of time
reflecting Brownian motion $X$ in $D$ spends in $B_k$ before
exiting the ball $B(z_k, 3c_0r_k/4)$ is bounded by $c_6 r_k^2$. By
the support theorem for standard Brownian motion, there exists
$p_2>0$ such that starting at any point
in
 $\partial B(z_k, 3c_0r_k/4)$, there is probability at least $p_2$
of hitting $\gamma_k$ before hitting
$B_k$.
 So the number of crossing from $B(z_k, 3c_0r_k/4)$ to $B_k$ by
reflecting Brownian motion $X^{(k)}$ in $\Omega_k$ killed upon
hitting $\gamma_k$ is majorized by a geometric random variable
with mean
$1/p_2$.
 This implies that the expected amount of time spent in
$B_k$
 by $X^{(k)}$ starting at any point
in
 $\gamma_{k+2}$ is at most $c_7 r_k^2$; that is
$$ \int_{B_k} G_{\Omega_k} (x, y) dx \leq c_7 r_k^2
\qquad \hbox{for } y\in  \gamma_{k+2}.
$$
Since
$|B_k|=c_{3}r_k^d$, the Harnack inequality implies that
$$ G_{\Omega_k}(z_k,y)\leq c_{8} r_k^{2-d}
\qquad \hbox{for every } y\in \gamma_{k+2}.
$$
Again by the Harnack and the Neumann boundary Harnack inequality, we have
$$
G_{\Omega_k} (x,y)\leq c_9 r_k^{2-d} \qquad \hbox{for every }
x\in \gamma_{k+1} \hbox{ and } y\in \gamma_{k+2}.
\eqno(4.4)
$$

For $k\geq 0$, it follows from the strong Markov property that
$$ G_{\Omega_k}(x, y)= G_{\Omega_{k+1}}(x, y)
+ \E_x \left[ G_{\Omega_k} (X_{T_{\gamma_{k+1}}}, y)
    \right] \qquad \hbox{for } x, y \in \Omega_{k+1}.
$$
Consequently, for every $x\in \Omega_{k+1}$ and $y\in \overline \Omega_{k+2}$.
$$ G_{\Omega_k}(x, y)= G_{\Omega_{k+1}}(x, y)
+ \E_x \left[ \E_y\left[ G_{\Omega_k} (X_{T_{\gamma_{k+1}}},
Y_{T^Y_{\gamma_{k+2}}}) \right] \right] , \eqno(4.5)
$$
where $Y$ is a reflecting Brownian motion in $D$ independent of
$X$ and $T^Y_{\gamma_{k+2}}$ is the first hitting time of
$\gamma_{k+2}$ by $Y$.
Let $a_k(x,y) = \E_x \left[ \E_y\left[ G_{\Omega_k}
(X_{T_{\gamma_{k+1}}}, Y_{T^Y_{\gamma_{k+2}}}) \right] \right]$
and note that by (4.3) and (4.4), for $x\in \Omega_{k+1}$ and
$y\in \overline \Omega_{k+2}$,
 $$c_5 r_{k}^{2-d} \leq a_k(x,y) \leq c_9 r_{k}^{2-d}.\eqno(4.6)$$
Fix $n\geq 1$. By (4.5) and (4.6), for every  $x\in \gamma_n$
and $y \in \gamma_{n+1}$ and $0\leq k\leq n-2$,
 $$G_{\Omega_k}(x, y) - G_{\Omega_{k+1}}(x, y) = a_k(x,y).$$
Adding these equations for $0\leq k\leq n-2$, we obtain
 $$G_{\Omega_0}(x, y) - G_{\Omega_{n-1}}(x, y)
 = \sum_{k=0}^{n-2} a_k(x,y),$$
or
  $$G(x,y) = G_{\Omega_0}(x, y)
  = G_{\Omega_{n-1}}(x, y) + \sum_{k=0}^{n-2} a_k(x,y).\eqno(4.7)$$
By (4.3) and (4.4), $c_5 r_{n-1}^{2-d} \leq G_{\Omega_{n-1}}(x, y)
\leq c_9 r_{n-1}^{2-d}$. This, (4.6) and (4.7) imply that
$$ c_5 \sum_{k=0}^{n-1} r_{k}^{2-d} \leq G(x, y) \leq c_9
\sum_{k=0}^{n-1} r_{k}^{2-d},
$$
for $x\in \gamma_n$ and $y \in \gamma_{n+1}$. By Definition 4.1,
$\alpha_1 <r_{n-1}/r_{n} < \alpha_2$ and $\alpha_1 <r_{0}/r_{1} <
\alpha_2$, so
$$ c_{10} \sum_{k=1}^{n} r_{k}^{2-d} \leq G(x, y) \leq c_9
\sum_{k=1}^{n} r_{k}^{2-d},
$$
for $x\in \gamma_n$ and $y \in \gamma_{n+1}$.
 \qed

\bigskip
\noindent{\bf Lemma 4.5}. {\sl For $n\geq 3$ and $y\in
\gamma_{n-3}$, let $x\mapsto h_y(x)$ be the Poisson kernel with
pole at $y$ for reflecting Brownian motion in $\Omega_{n-3}$
killed upon hitting $\gamma_{n-3}$. Let $\P_x^y$ denote the
distribution of the $h_y$-transform of reflecting Brownian motion
in $D^*$ killed upon hitting $\gamma_{n-3}$, starting from $x\in
\Omega_{n-3}$. There exist $c_1,p_1>0$, depending only on $D$,
such that $\P_x^y(L_{T_{\gamma_{n-3}}} > c_1 r_n)>p_1$ for any
collection of $\gamma_k$'s as in Definition 4.1, any $n\geq 3$,
$x\in \gamma_n$ and $y\in \gamma_{n-3}$. }

\bigskip
\noindent{\bf Proof}. All constants $c_j$ that appear in this
proof depend only on $D$. The conditions listed in Definition 4.1
imply existence of $c_2>0$ and a point $x_0\in \prt D_{n-2} \cap
\prt D$ such that (i) the distance from $x_0$ to $\gamma_{n-1}
\cup \gamma_{n-2}$ is greater than $2c_2 r_n$, and (ii) there
exists an orthonormal coordinate system $CS_{x_0}$ such that
$B(x_0, c_2 r_n)\cap \prt D$ is the graph of a Lipschitz function
with the Lipschitz constant $\lambda$. Recall that $\lambda$ is
the constant in the definition of $\cD$ and, hence, in the
definition of $\cD_1$. We will assume that $x_0 = 0$ in $CS_{x_0}$
and the positive part of the $d$-th coordinate axis intersects $B(x_0, c_2
r_n)\cap D$.

Let $h(x) = \P_x(T_{\gamma_{n-3}} < T_{\gamma_{n}})$ and let $x_1$
be the intersection point of $ \prt B(x_0, c_2 r_n/2)$ and the
positive part of the $d$-th coordinate axis in $CS_{x_0}$. It is
easy to show, using Definition 4.1 and a ``Harnack chain of
balls'' argument, that $h(x_1)> c_3>0$. By the boundary Harnack
principle (Lemma 2.8), we have $h(x) > c_4>0$ for all $x\in B(x_0,
3c_2 r_n/4) \cap D$.

Recall that $X^*$ is reflecting Brownian motion on the
Martin-Kuramochi compactification $D^*$ of $D$, and $X$ is the
quasi-continuous projection of $X^*$ into $\overline D$. As we
noted in Section 2, $X=X^*$ on $D\cup \prt_L D$. Let $X^*$
start from a point $x\in B(x_0, c_2 r_n / 2)\cap \prt D \subset
\prt_L D$. It follows from Theorem 2.2 that $X_t = x + W_t + N_t$,
where $W_t$ is a $d$-dimensional Brownian motion starting from 0
and $N_t=\int_0^t {\bf n} (X_s) dL_s$ is the singular push on the
boundary $\prt_L D$. Assume without loss of generality that
$\lambda>1$. Let $x_2=(0,0,\dots, 0, -c_2 r_n /10)$, $B_1 = B(x_2,
c_2 r_n / (100\lambda))$, $B_2=B(0, c_2 r_n /(100\lambda))$, and
let $C_1$ be the convex hull of $B_1 \cup B_2$. Consider the event
$A$ that the Brownian motion $W$ hits $B_1$ before leaving $C_1$
in less than $c_5 r_n^2$ units of time. By the support theorem and
Brownian scaling, the probability of $A$ is greater than $p_2>0$.
Let $T_* = T^X_{\prt B(x_0, 3c_2 r_n/4)} \land c_5 r_n^2$. We will
argue that if $A$ occurs, then $|N_{T_*}| \geq c_6 r_n$. To see
this, first suppose that $T^X_{\prt B(x_0, 3c_2 r_n/4)} < c_5
r_n^2$. Since $A$ holds, $W$ stays in $C_1$, and it follows that
$|W_{T_*}| \leq c_2 r_n /5$. Since $|X_{T_*} - x| \geq c_2 r_n/4$,
we have $|N_{T_*}| \geq c_2 r_n/20$. If $T^X_{\prt B(x_0, 3c_2
r_n/4)} \geq c_5 r_n^2$ and $A$ holds, let $t_0\leq T_*$ be a time
when $W_{t_0} \in B_1$. Since $x+ B_1$ is at a distance greater
than $c_2 r_n/(100\lambda)$ from $D$ and $X_{t_0} \in \ol D$, we
must have $|N_{T_*}| \geq c_2 r_n/(100\lambda)$. We see that with
probability $p_2$ or greater, $X$ accumulates at least $c_6 r_n$
units of local time before leaving $B(x_0, 3c_2 r_n/4)$. Since
$h(x) > c_4$ for all $x\in B(x_0, c_2 r_n) \cap D$, $X$ starting
from any point $x\in B(x_0, c_2 r_n / 2)\cap \prt D$ has a chance
$p_3>0$ (depending only on $D$) of accumulating at least $c_6 r_n$
units of local time and hitting $\gamma_{n-2}$ before hitting
$\gamma_{n}$, by the strong Markov property applied at the time of
hitting of $\prt B(x_0, 3c_2 r_n/4)$.

Suppose that $x_1 \in \gamma_{n-1}$ lies at least $c_7 r_n$ units
away from $\prt D$. By the support theorem for Brownian motion,
the chance that reflecting Brownian motion $X$ starting from
$x_1$ will hit $B(x_0, c_2 r_n/2) \cap \prt D$ before hitting any
other part of $\prt D \cup \gamma_n \cup \gamma_{n-2}$ is greater
than $p_4>0$, depending only on $D$. By the strong Markov property
applied at the hitting time of $B(x_0, c_2 r_n) \cap \prt D$,
reflecting Brownian motion starting from $x_1$ has a chance
greater than $p_5>0$ of accumulating at least $c_6 r_n$ units of
local time inside $B(x_0, 3c_2 r_n/4)$ and hitting $\gamma_{n-2}$
before hitting $\gamma_{n}$.
Hence, the $h$-transform of $X$
 starting from $x_1$ has a chance greater than $p_5>0$ of
accumulating at least $c_6 r_n$ units of local time inside $B(x_0,
3c_2 r_n/4)$ before its lifetime. The boundary Harnack principle
shows that the same is true for any $x\in \gamma_{n-1}$, except
that the bound for the probability has to be replaced with a new
value $p_6>0$.

Now consider
the $h_y$-transform of $X$
 starting from a point of $\gamma_n$. It must hit $\gamma_{n-1}$
and then $\gamma_{n-2}$ on its way to $\gamma_{n-3}$. The strong
Markov property applied at the hitting times of $\gamma_{n-1}$ and
$\gamma_{n-2}$ and the claims proved so far show that an
$h$-process starting from any point in $\gamma_n$ has a chance
$p_6>0$ of accumulating at least $c_6 r_n$ units of local time
inside $B(x_0, 3c_2 r_n/4)$ before its lifetime. By the Neumann
boundary Harnack principle
(Lemma 2.8), there are positive constants $c_7<c_8$ such that
$$
 c_7 \frac{h(x)}{h(z)}\leq  \frac{h_y(x)}{h_y(z)} \leq c_8
 \frac{h(x)}{h(z)} \qquad
\hbox{for } x, z\in \gamma_{n-2} \cup \gamma_{n-1}.
$$
A routine argument based on this observation allows us to extend the
claim to the $h_y$-process.
 \qed

\bigskip
\noindent{\bf Proof of Theorem 4.3}. The main idea of this proof
is the same as that of the proof of Theorem 3.2 but some details
are different.

(i) Consider $x_0 \in D$. It is not hard to see that there exists
$z_0\in \prt D$ and a corresponding family of $\gamma_n$'s, as in
Definition 4.1, such that $x_0 \in D_{n_0+1}$ for some $n_0\geq 1$
and
$\dist(x_0, \prt D) \geq c_1 r_{n_0+1}$.
 For $y\in \Omega'_{n_0-1}$, by the strong Markov property of $X$,
$$ G(x_0, y) =\E_{x_0} \left[ G(X_{T_{\gamma_{n_0}}}, y) \right]
=\E_{x_0} \left[ \E_y \left[ G\Big(X_{T_{\gamma_{n_0}}},
Y_{T^Y_{\gamma_{n_0-1}}} \Big);
  \, T^Y_{\gamma_{n_0-1}}< T^Y_{B_*}   \right] \right],
$$
where $Y$ is a reflecting Brownian motion in $D^*$ killed upon hitting
$\partial B_*$ starting from $y$ and
independent of $X$, and $T^Y_{\gamma_k}$ is the first hitting
time
 of $\gamma_{k}$ by $Y$. Hence by Lemma 4.4, $G(x_0, y)$ is bounded
by $c_2 \sum_{k=1}^{n_0} r_k^{2-d}$ for $y\in \Omega'_{n_0-1}$. By
the Harnack principle, $y\mapsto G(x_0, y)$ is bounded by $c_3
\sum_{k=1}^{n_0} r_k^{2-d}$ on
$\prt B(x_0, c_1 r_{n_0+1}/2)$,
 and the maximum principle implies that the same bound holds on
$D \setminus B(x_0, c_1 r_{n_0+1}/2)$.
 We obtain, with $\sigma$ denoting the surface measure on $\prt_L
D=\bigcup_{n=1}^\infty ( \partial D_n \cap \partial D) $,
$$\eqalign{
 \E_{x_0} L_{T_{B_*}}
&=\E_{x_0} \left[ \int_0^{T_{B_*}} \bone_{\prt_L D}(X_t) dL_t \right]\cr
&=\int_{\prt_L D} G(x_0, x) \sigma (dx) \cr
&= \sum_{n=1}^{n_0 -1}
 \int_{\prt D_n \cap \prt D} G(x_0, x) \sigma (dx)
 +\sum_{n= n_0}^\infty  \int_{\prt D_n \cap \prt D}
 G(x_0, x) \sigma (dx)\cr
 &\leq  \sum_{n=1}^{ n_0-1} |\prt D_n \cap \prt D|
 \sup_{x\in \prt D_n \cap \prt D} G(x_0,x)
 + \sum_{n=n_0}^\infty |\prt D_n \cap \prt D|
 \sup_{x\in \prt D_n \cap \prt D} G(x_0,x)\cr
 &\leq  \sum_{ n=1}^{n_0-1} |\prt D_n \cap \prt D|
\,  c_4 \sum_{k=1}^{n} r_k^{2-d}
 + \sum_{n= n_0}^\infty |\prt D_n \cap \prt D|
 \, c_3 \sum_{k=1}^{n_0} r_k^{2-d}\cr
 &\leq  \sum_{n=1}^\infty  c_5 |\prt D_n \cap \prt D|
 \sum_{k=1}^{n} r_k^{2-d}.
 }$$
This is bounded by a constant independent of $x_0$, by assumption
(4.1). Hence we obtain $\sup_{x\in D} \E_{x} L_{T_{B_*}} < \infty$
and, therefore, $ \inf_{x\in D} u(x) =\inf_{x\in D} \E_x
\exp(-L_{T_{B_*}}) >0$. This means that the whole surface of $D$
is active.

(ii) Consider a point $z_0 \in \prt D$ and a corresponding family
of $\gamma_n$'s satisfying (4.2). Let $\{x_n, n\geq 1\}$ be a
sequence in $D$ that converges to $z_0$. There is a subsequence
$\{x_{n_j}, j\geq 1\}$ that converges to some $z_0^*$ in $D^*$,
the Martin-Kuramochi compactification of $D$. Note that $x\mapsto
G(z_0^*, x)$ is well defined on $D^*\setminus \{z_0^*\}$ by the
second paragraph of Section 3. In particular,
$$G(z_0^*, x)=\lim_{j\to \infty} G(x_{n_j}, x) \qquad
\hbox{for } x\in D^*\setminus \{z_0^*\}.
$$

For $a\geq 0$, define
$$\eta_a =\{x\in D^*\setminus (B_*\cup \{z_0^*\}): G(z^*_0, x) = a\}.
$$
Note that $\eta_0=\prt B_*$. Recall the definition of the
integer $m_0 > 1$ from
Condition  4.2. Define for $n\geq 1$,
$$ a_n = \inf_{x\in \gamma_{3nm_0}} G(z_0^*, x).
$$
By the Neumann boundary Harnack principle applied to the function
$x\mapsto G(z_0^*, x)$ on $\gamma_{3nm_0}$, there is a constant
$c_0\in (1,\infty)$ depending only on $D$ such that $\inf_{x\in
\gamma_{3nm_0}} G(z_0^*, x) \geq c_0 \E_x [ G(z_0^*,
X_{T_{\gamma_{3nm_0}}} )]$ for every $x\in \gamma_{3(n+1)m_0}$.
For $n\geq 0$,
$$\eqalignno{
a_{n+1}-a_n &= \inf_{x\in \gamma_{3(n+1)m_0}} G(z_0^*, x)
   - \inf_{x\in \gamma_{3nm_0}} G(z_0^*, x) \cr
&\leq \inf_{x\in \gamma_{3(n+1)m_0}} \left(  G(z_0^*, x)
   - c_0 \E_x \left[ G(z_0^*, X_{T_{\gamma_{3nm_0}}} ) \right]
   \right)\cr
&\leq \inf_{x\in \gamma_{3(n+1)m_0}} \left( c_0 G(z_0^*, x)
   - c_0 \E_x \left[ G(z_0^*, X_{T_{\gamma_{3nm_0}}} ) \right]
   \right)\cr
 &=c_0 \inf_{x\in \gamma_{3(n+1)m_0}} G_{\Omega_{3nm_0}} (z_0^*, x)
 \cr
 &=c_0 \inf_{x\in \gamma_{3(n+1)m_0}} \E_{z_0^*} \left[
 G_{\Omega_{3nm_0}} (X_{T_{\gamma_{3(n+2)m_0}}}, x)  \right] \cr
 &\leq c_1 r_{3nm_0}^{2-d}, &(4.8) \cr}
$$
where the last inequality is due to
(4.4)
 and $c_1>0$ is a constant depending only on $D$.

On the other hand, for $x\in \gamma_n$,
$$
G(z_0^*, x)= \lim_{j\to \infty} G(x_{n_j}, x)
= \lim_{j\to \infty} \E_{x_{n_j}}
\left[ G\big(X_{T_{\gamma_{3(n+1)m_0}}} , x \big) \right].
$$
So it follows from Lemma 4.4 that  there are positive constants
$c_2<c_3$ depending only on
$D$
 such that
$$ c_2 \sum_{k=1}^{3nm_0} r_k^{2-d} \leq a_n=\inf_{x\in \gamma_{3nm_0}}
G(z_0^*, x) \leq c_3 \sum_{k=1}^{3nm_0} r_k^{2-d}
 \qquad \hbox{for every } n\geq 1. \eqno (4.9)
$$
It follows from Definition 4.1 that for some $c_4 < \infty$,
$r_n^{2-d}/r_{n-1}^{2-d} \leq c_4$ for every $n\geq 1$. Hence,
 $$\eqalignno{{a_n\over a_{n-1}} &\leq {
 c_3
 \sum_{k=1}^{3nm_0} r_k^{2-d} \over
 c_2 \sum_{k=1}^{3(n-1)m_0} r_k^{2-d}}
 \leq {c_3\over c_2} + {c_3 \sum_{k=3(n-1)m_0+1}^{3nm_0} r_k^{2-d} \over
 c_2  r_{3(n-1)m_0}^{2-d}} &(4.10)\cr
& \leq {3c_3m_0\over c_2} ( 1 + c_4^{3m_0})  <
 \infty. \cr}
 $$
Let $\beta \df \max \{ 2, \, (3c_3m_0/ c_2) ( 1 + c_4^{3m_0})\}$
and $n_j =\inf\{n: a_n \geq \beta^{j}\}$. By (4.10),
$$ \beta^j\leq a_{n_j} \leq a_{n_j-1} \beta \leq \beta ^{j+1}.
\eqno(4.11)
$$
Since by (4.2),
$$
\sum_{m=0}^{3m_0 -1} \sum_{i=0}^1 \sum_{j=1}^\infty \ \sum_{\{n:
3nm_0+m \in (n_{4j-4+2i}, n_{4j-2+2i}]\}}
r_{3nm_0+m}^{d-1}
 \sum_{k=1}^{3nm_0+m}
  r_k^{2-d}=\infty,
$$
without loss of generality, we may and do assume the sum is
infinite for $m=0$ and $i=1$, i.e.,
 $$
\sum_{j=1}^\infty \ \sum_{\{n: 3nm_0 \in (n_{4j-2}, n_{4j}]\}}
r_{3nm_0}^{d-1}
 \sum_{k=1}^{3nm_0}
  r_k^{2-d}=\infty. \eqno(4.12)
$$
Let  $X^*$ be reflecting Brownian motion on $D^*$ starting from
some $x_0\in D_{n_0}$, where $n_0$ is large. Define
 $$\eqalign{
 S_j &=\inf\{t>0: X^*_t \in \eta_{a_{n_j}}\}, \quad j\geq 1,\cr
T^{j,n}_1 &= \inf\{t>S_j: X^*_t \in\eta_{a_{n}}\}, \quad n\geq 1,
\cr U^{j,n}_k &= \inf\{t>T^{j,n}_k: X^*_t \in\eta_{a_{n-1}}\},
\quad
 n,k\geq 1, \cr
 T^{j,n}_k &= \inf\{t>U^{j,n}_{k-1}: X^*_t \in\eta_{a_{n}}\}, \quad
 n,k\geq 2, \cr
 N^j_n &=\max \{ k: U^{j,n}_k \leq S_{j-1}\},\quad n\geq 1.
 }$$
In words, $N^j_n$ is the number of downcrossings of
$[a_{n-1},a_{n}]$ by $M_t\df G(z_0^*, X^*_t)$ between times $S_j$
and $S_{j-1}$. This is of interest to us only for $n$ such that
$[a_{n-1},a_{n}] \subset [a_{n_{j-1}}, a_{n_j}]$.

The process $M$ is a continuous local martingale so it is a time-change of
Brownian motion, until it hits $0$.

Consider a one-dimensional Brownian motion $W_t$ starting from
$a_{n_{4j}}$ and killed at the hitting time $T$ of $a_{n_{4j-4}}$.
Note that by (4.11), $a_{n_{4j}}\geq \beta^{4j}$, $a_{n_{4j-4}}
\leq \beta^{4j-3}$, and $a_{n_{4j-2}} \geq \beta^{4j-2}$. It
follows from the Ray-Knight theorem that there is an event $A$
with probability greater than $p_1>0$, such that on $A$, the local
time $L^x_T$ accumulated by $W$ at level $x$ before time $T$ is
greater than $c_5 \beta^j$ for all $x\in (a_{n_{4j-2}} ,
a_{n_{4j}})$. We will apply excursion theory to excursions of $W$
from the set $\{a_{n}:\,  n_{4j-2} \leq n \leq n_{4j}\}$. Given
the local time $\{L^x_T: \,
x=a_{n}\in \{a_{n_{4j-2}} , a_{n_{4j}}\} \}$
 and assuming the event $A$ occurs, the distribution of the number
of excursions going from $a_{n-1}$ to $ a_{n}$ is minorized by a
Poisson random variable with expectation $K_j\geq c_5 {\beta^{4j}
\over a_{n}-
a_{n-1} }$.
 By (4.8), we have
$$ K_j \geq {c_5\over c_1} {\beta^{4j}\over r_{3nm_0}^{2-d}}.
$$
Conditional on $\{L^x_T, x=a_{n} \in [a_{n_{4j-2}} ,
a_{n_{4j}}]\}$, these random variables are independent.

Let $T^M(b) = \inf\{t>0: M_t = b\}$ and $M^j_t = \{M_{t+
T^M(a_{n_{4j}})}, t\in [0, T^M(a_{n_{4j-4}}) -
T^M(a_{n_{4j}})]\}$. Since $M^j_t$ is a time-change of $W_t$,
there exists an event $A'$ with $\P_{x_0}(A')>p_1$, such that on
$A'$, conditional on the local time of $M^j$
at $\{a_{n}, n_{4j-1} \leq n \leq n_{4j}\}$,
the numbers of excursions of $M^j$ between consecutive points of
$\{a_{n}, n_{4j-1} \leq n \leq n_{4j}\}$ are independent random
variables minorized by independent Poisson random variables with
means $K_n\geq {c_5\over c_1} {\beta^{4j}\over r_{3nm_0}^{2-d}}$.
By (4.9) and (4.11), there is a constant $c_6>0$ depending only
on $D$ such that
$$ K_j \geq c_6 {\sum_{k=1}^{n_{4j}} r_k^{2-d} \over r_{3nm_0}^{2-d}}
$$

Note that the processes $M^j$ are independent.
We will now consider the process $X$ conditioned on the values of
the following random elements: local times of the $M^j$'s
accumulated at levels $\{a_{n}, n_{4j-1} \leq n \leq n_{4j}\}$,
and the endpoints of excursions of $X$ from $\{\eta_{a_{n}},
n_{4j-1} \leq n \leq n_{4j}\}$.

It follows from the definition of $a_n$ and Condition  4.2
that
the following condition is satisfied

$$ \hbox{\sl There are at least three consecutive  $D_k$'s
between $\eta_{a_{n-1}}$ and $ \eta_{a_n}$.} \eqno (4.13)$$
\medskip

This and an easy argument based on Lemma 4.5 show that given
endpoints of an excursion of $X$ going from $\eta_{a_n}$ to
$\eta_{a_{n-1}}$, the amount of local time accumulated by the
excursion on $\prt D$ is greater than $c_7 r_{3nm_0}$ with
probability greater than $p_2>0$.

Let $J_n$ be the distribution of the local time accumulated by $X$
on the part of $\prt D$ between $\eta_{a_{n-1}}$ and $\eta_{a_n}$
during the time interval $(T^M(a_{n_{4j}}), T^M(a_{n_{4j-4}}))$.
We have shown that on an event $A_j$ of probability greater than
$p_1$, $J_n$ is stochastically minorized by a random variable
$I_n$ whose distribution is Poisson with mean greater than
$$
 p_2 \, c_6 \, {\sum_{k=1}^{n_{4j}} r_k^{2-d} \over r^{2-d}_{3n
 m_0}} \,
 c_7 r_{3nm_0} = p_2\, c_6 \, c_7\, r_{3nm_0}^{d-1}
 \sum_{k=1}^{n_{4j}} r_k^{2-d}.
$$
Hence $J_n$ is minorized by a random variable with mean
$\lambda_n\geq c_8 r_{3nm_0}^{d-1} \sum_{k=1}^{n_{4j}} r_k^{2-d}$
and variance $\lambda_n$. Moreover, we can assume that the $I_n$'s
are independent given $A_j$. Hence, the local time accumulated by
$X$ between the hitting of $\eta_{a_{n_{4j}}}$ and
$\eta_{a_{n_{4j-2}}}$, on the part of $\prt D$ between these
surfaces, is stochastically minorized by a random variable $H_j$
such that on the event $A_j$, its mean is bounded below by
$$
 \sum_{\{n: 3nm_0 \in [n_{4j-2}, n_{4j}]\}} c_8 r_{3nm_0}^{d-1}
 \sum_{k=1}^{n_{4j}} r_k^{2-d}
$$
and the variance equals its mean. It follows that $H_j$ takes a
value no less than
$$ b_j\df {1\over 2} \sum_{\{n: 3nm_0 \in [n_{4j-2}, n_{4j}]\}} c_8
 r_{3nm_0}^{d-1}
  \sum_{k=1}^{n_{4j}} r_k^{2-d}
$$
with probability greater than $p_3>0$.

Since the $M^j$'s are independent, we can assume that the $H_j$'s
are independent. Let $\Lambda_j$ be independent random variables
with $P(\Lambda_j =b_j) = 1- P(\Lambda_j = 0) = p_3$. The
distribution of the local time accumulated by reflecting Brownian
motion starting from $x_0\in D_{n_{4j_0}}$ before hitting $B_*$ is
minorized by the distribution of $\sum_{j=1}^{j_0} \Lambda_j$. In
view of (4.12), $\sum_{j\geq 1} b_j = \infty$, and this easily
implies that $\sum_{j \geq 1}\Lambda_j=\infty$, a.s. Hence, for
any $b< \infty$, there is some  $j_0$ such that $\P (\sum_{j \leq
j_0}\Lambda_j>b)> 1-1/b$. This implies that $\P_{x_0}
(L_{T_{B_*}}>b)> 1-1/b$ for every $x_0\in D_{n_{4j_0}}$.
Therefore, $\inf_{x\in D}\E_x\left[ \exp(-L_{T_{B_*}})\right]=0 $
and we see that part of the surface of $D$ is nearly inactive.
 \qed

\bigskip

We note that the only place where Condition  4.2 is used in
the proof of Theorem 4.3(ii) is to prove (4.13).
 We will next discuss Condition 4.2 but first we need a
lemma.

\bigskip

\noindent{\bf Lemma 4.6}. {\sl Let $D\in \cD_1$, $z\in \prt D$
and let $\{D_n, n\geq 0\}$ and  $\{\gamma_n, n\geq 0\}$ be as in
 Definition 4.1. There exist $\alpha_0,p_0>0$ depending only on $D$
such that the following holds.

\item{(i)} For every $n\geq 1$ and every positive harmonic
function $h$ on the interior of $\ol{D_{n-1} \cup D_n}$, with
Neumann boundary conditions on $\prt D \cap \ol{D_{n-1} \cup
D_n}$,
$$ h(x) \leq \alpha_0 h(y) \qquad \hbox{for every } x,y\in \gamma_n.
$$

\item{(ii)} For every $n\geq 1$,
$$  \P_x (T_{\gamma_{n+1}} <
T_{\gamma_{n-1}}) \geq p_0 \qquad \hbox{for every } x \in
\gamma_n.
$$
}

\pf
(i)
 Let $h$ be a positive harmonic function on the interior of
$\ol{D_{n-1} \cup D_n}$, with Neumann boundary conditions on $\prt
D \cap \ol{D_{n-1} \cup D_n}$.
It follows from parts (i) and (iv) of Definition 4.1 that there
are $k_1<\infty$ and $c_1>0$, depending only on $D$, such that
$\gamma_n$ can be covered by at most $k_1$ balls $B(x_k, c_1
r_n)$. Moreover, for each one of these balls, either $B(x_k, 2c_1
r_n)\subset D$ or $B(x_k, 2c_1 r_n)\cap \prt D$ is the graph of a
Lipschitz function.
If $x,y\in \gamma_n$ and both points belong to one of balls
$B(x_k, c_1 r_n)$,
 then there is a constant
$c_2>0$
 depending only on $D$ such that
$h(x) \leq c_2 h(y)$,
 either by the usual Harnack principle (if the ball is inside $D$)
or by the Neumann boundary Harnack principle proved in Lemma 2.8.
It follows by a Harnack chain argument that $h(x) \leq \alpha_0
h(y)$ for any $x,y\in \gamma_n$, where
$\alpha_0= c_2^{k_1}$.

(ii)
 According to the definition of $D\in \cD_1$,
 there is $k_2<\infty$ such that there exists a ``Harnack chain of
balls'' connecting
 $\gamma_n$ and $\gamma_{n+1}$, that is,
we can find
 a sequence $ B(x_1,r), B(x_2,r), \dots, B(x_k, r)$ in
$\Omega_{n-1}$ with
$k\leq k_2$,
 $x_1 \in \gamma_n$, $x_k \in \gamma_{n+1}$ and $x_j \in B(x_{j-1},
r/2)$ for $j=2,\dots , k$. The existence of this ``Harnack chain
of balls''  and the Harnack inequality easily imply that for some
$p_1>0$ depending only on $D$ and some $x \in \gamma_n$,
$$ \P_x (T_{\gamma_{n+1}} < T_{\gamma_{n-1}}) \geq p_1.
$$
Applying part (i) of this
lemma
 to harmonic function $x\mapsto \P_x (T_{\gamma_{n+1}} <
T_{\gamma_{n-1}})$, we conclude that there is some $p_0>0$
depending only on $D$ such that $\P_x (T_{\gamma_{n+1}} <
T_{\gamma_{n-1}}) \geq p_0 $ for every  $x\in \gamma_n$. \qed

\bigskip

Part (ii) of Theorem 4.3 has been proved under the assumption that
Condition  4.2 holds. Condition  4.2 seems to be difficult to
verify in a direct way. We will state two other conditions,
Conditions 4.7 and  4.8, that are easier to verify in examples. We
will show that Condition 4.7 implies Condition 4.8 and Condition 4.8
implies Condition 4.2. In some examples,  Condition 4.7 is the easiest
condition to verify, but in some other examples  Condition 4.8 holds
even though  Condition 4.7 does not. Lemma 4.11 below shows how one
can verify Condition  4.7 in some examples.

\bigskip

\noindent {\bf  Condition 4.7.}
{\sl Let $\alpha_0$ and $p_0$ be the constants in Lemma 4.6 and
let $T_{\gamma_n}$ be the first hitting time of $\gamma_n$ by
reflecting Brownian motion in $D$.
 There exist $0<m_0<m_1\leq \infty$ such that for any $z\in\prt D$
and the $\gamma_n$'s as in Definition 4.1 corresponding to $z$, if
$n>m_1$,
$$\P_x (T_{\gamma_{n-m_0-1}} < T_{\gamma_{n+1}}) \leq \alpha_0^{-2} p_0
\qquad \hbox{for every } x \in \gamma_n,
$$
and
$$ \P_x (T_{\gamma_{n+1}} < T_{\gamma_{n-m_0-1}}) \leq \alpha_0^{-2} p_0
\qquad \hbox{for every } x\in \gamma_{n-m_0}.
$$
}

\bigskip

\noindent {\bf  Condition 4.8.} {\sl There exist $0<m_0\leq m_1<\infty$ such that
for any $z\in\prt D$ and the $\gamma_n$'s as in Definition 4.1, the
following is true for $n>m_1$. Let $A$ be the interior of
$\ol{\bigcup_{n-m_0-1 \leq k \leq n} D_k}$ and
$\mu_x(dy)=\P_x(T_{\gamma_{n-m_0-1}\cup\gamma_{n+1}} \in dy)$, for
$x\in A$. In other words, $\mu_x$ is harmonic measure on the set
$\gamma_{n-m_0-1}\cup\gamma_{n+1}$ inside $A$ for Brownian motion
reflected on $\prt D$. Then the Radon-Nikodym derivative
$d\mu_z/d\mu_y \leq 1$ on $\gamma_{n+1}$ and $d\mu_z/d\mu_y \geq
1$ on $\gamma_{n-m_0-1}$, for $z \in \gamma_{n-m_0}$ and  $y \in
\gamma_{n}$.
}

\bigskip

\noindent{\bf Lemma 4.9}. {\sl Condition  4.7 implies
Condition 4.8.
}

\bigskip

\noindent{\bf Proof}. Recall
 the constants $\alpha_0$ and $p_0$ from
Lemma
 4.6.

Let $A_n $ be the interior of $\ol{\bigcup_{n-m_0-1 \leq k \leq n}
D_k}$ and $\mu_x(dy)=\P_x(X_{T_{\gamma_{n-m_0-1}\cup\gamma_{n+1}}}
\in dy)$, for $x\in A_n$. In other words, $\mu_x$ is harmonic
measure in $A_n$ for Brownian motion reflected on $\prt D$. Fix a
set $C\subset \gamma_{n+1}$. By Lemma 2.7, $x\mapsto \mu_x(C)$ is
a non-negative harmonic function of $x\in A_n$. By
Lemma
 4.6,
$$\mu_x(C)\leq \alpha_0 \mu_y(C) \qquad \hbox{and} \qquad
\mu_y(\gamma_{n+1}) \leq
\alpha_0
 \mu_x(\gamma_{n+1})
$$
for $x,y\in \gamma_n$, and so
$$\mu_x(C) \leq \alpha_0^2 \,
{ \mu_y(C) \over\mu_y(\gamma_{n+1})} \, \mu_x(\gamma_{n+1}) \qquad
\hbox{for } x,y\in \gamma_n. \eqno(4.14)
$$
By
Lemma
 4.6 and Condition 4.7, for $n>m_1$,
$\mu_x(\gamma_{n+1}) \geq p_0$ for all $x \in \gamma_{n}$ and
$\mu_z(\gamma_{n+1}) \leq \alpha_0^{-2}p_0$ for all $z \in
\gamma_{n-m_0}$. Hence,
$$
\mu_z(\gamma_{n+1}) /\mu_x(\gamma_{n+1}) \leq \alpha_0^{-2},
 \eqno(4.15)
$$
for all $z \in \gamma_{n-m_0}$, $x \in \gamma_{n}$ and $n>m_1$. If
reflecting Brownian motion in $D$ starts from a point in
$\gamma_{n-m_0}$, it has to hit $\gamma_n$ before hitting
$\gamma_{n+1}$. Hence, by the strong Markov property, (4.14), and
(4.15), we obtain for $n>m_1$, $C\subset \gamma_{n+1}$, $z \in
\gamma_{n-m_0}$ and $y \in \gamma_{n}$,
 $$\eqalign{
 \mu_z(C) &=
 \int_{\gamma_n} \mu_x(C) \P_z(
 X_{T_{\gamma_{n-m_0-1}\cup \gamma_n}}
 \in dx)\cr
 &\leq \alpha_0^2 \, {\mu_y(C) \over \mu_y(\gamma_{n+1})} \,
 \int_{\gamma_n} \mu_x(\gamma_{n+1})
 \P_z(
 X_{T_{\gamma_{n-m_0-1}\cup\gamma_n}}
 \in dx)\cr
 &= \alpha_0^2 \, {\mu_y(C) \over \mu_y(\gamma_{n+1})}
 \mu_z(\gamma_{n+1}) \cr
&\leq \mu_y(C) .
 }$$
Since $C$ is an arbitrary subset of $\gamma_{n+1}$, the
Radon-Nikodym derivative $d\mu_z/d\mu_y \leq 1$ on $\gamma_{n+1}$
for $z \in \gamma_{n-m_0}$ and  $y \in \gamma_{n}$ with $n>m_1$.
Similarly, $d\mu_z/d\mu_y \geq 1$ on $\gamma_{n-m_0-1}$ for $z \in
\gamma_{n-m_0}$ and  $y \in \gamma_{n}$ with $n>m_1$. Therefore
Condition  4.8 is satisfied.
 \qed

\bigskip

\noindent{\bf Lemma 4.10}. {\sl Condition  4.8 implies Condition
4.2. }

\bigskip

\noindent{\bf Proof}. We will consider $n\geq m_1+1$. Suppose that
$x_0 \in \Omega_{n+1}$ and choose $c_1$ so large that $K\df \{x\in
D\setminus B_*: G(x_0,x) \geq c_1\} \subset \Omega_{n+1}$. We will
prove that $\sup_{x\in \gamma_{n-m_0}} G(x_0,x) \leq \inf _{x\in
\gamma_n} G(x_0,x)$. It will suffice to show that $\P_x(T_K <
T_{B_*}) \leq \P_y(T_K < T_{B_*})$ for $x \in \gamma_{n-m_0}$ and
$y \in \gamma_{n}$, because $G(x_0,x) = c_1 \P_x(T_K < T_{B_*})$
for $x\in D\setminus(B_*\cup K)$.

Our proof will use the technique of coupling. We will construct
two reflecting Brownian motions in $D$ on a common probability
space, $X$ starting from $x \in \gamma_{n-m_0}$ and $Y$ starting
from $y \in \gamma_{n}$, such that $\{T^X_K < T^X_{B_*}\} \subset
\{T^Y_K<T^Y_{B_*}\}$ almost surely, that is,
$$
  \P \left(T^X_K < T^X_{B_*} \hbox{ and }  T^Y_K\geq T^Y_{B_*} \right) =0.
$$

Let $A_n$ be the interior of $\ol{\bigcup_{n-m_0-1 \leq k \leq n}
D_k}$ and define
$$\mu_x(dy)\df \P_x(X_{T_{\gamma_{n-m_0-1}\cup\gamma_{n+1}}} \in
dy) \qquad \hbox{for } x \in A_n.
$$
Given $x\in \gamma_{n-m_0}$ and $y \in \gamma_{n}$, we will define
some random variables on a common probability space. Let
$\eta_{x,n+1}$ be a random variable taking values in
$\gamma_{n+1}$ and having distribution
$\mu_x(dz)/\mu_x(\gamma_{n+1})$. Let $\eta_{y,x,n+1}$ take values
in $\gamma_{n+1}$ with distribution $(\mu_y(dz) - \mu_x(dz))/
(\mu_y(\gamma_{n+1}) -\mu_x(\gamma_{n+1}))$, and let $I_{x,n+1}$
take values 0 or 1, with $P(I_{x,n+1} = 1) = \mu_x(\gamma_{n+1})$.
Note that $\eta_{y,x,n+1}$ is well defined because, under
Condition  4.8, $\mu_y(dz) \geq \mu_x(dz)$ on
$\gamma_{n+1}$. Similarly, let $\eta_{y,n-m_0-1}$ be a random
variable taking values in $\gamma_{n-m_0-1}$ and having
distribution $\mu_y(dz)/\mu_y(\gamma_{n-m_0-1})$ on
$\gamma_{n-m_0-1}$. Let $\eta_{x,y,n-m_0-1}$ take values in
$\gamma_{n-m_0-1}$ and have distribution $(\mu_x(dz) -
\mu_y(dz))/(\mu_x(\gamma_{n-m_0-1})-\mu_y(\gamma_{n-m_0-1}))$. Let
$I_{y,n-m_0-1}$ take values 0 or 1, and assume that
$P(I_{y,n-m_0-1} = 1) = \mu_y(\gamma_{n-m_0-1})$. Due to Condition
4.8, we may and do assume that the $I$'s are constructed so
that $I_{x,n+1} + I_{y,n-m_0-1} \leq 1$, a.s., and we let $I_{x,y}
= 1 - I_{x,n+1} - I_{y,n-m_0-1}$. Moreover, the $\eta$'s are
constructed so that they are independent, and independent of
the $I$'s.

For $x\in A_n$ and $z\in \gamma_{n-m_0-1}\cup\gamma_{n+1}$, let
$\Q_x^z$ denote the distribution of reflecting Brownian motion $X$
in $A_n$ starting from $x$, conditioned on leaving $\overline
{A_n}\setminus\{ \gamma_{n-m_0-1}\cup\gamma_{n+1}\}$ through $z$.
Let $\Q_{x,y}^z$ denote the distribution of a pair of processes
$(\wh X,\wh Y)$, such that the distribution of $\wh X$ is $\Q_x^z$
and the distribution of $\wh Y$ is $\Q_y^z$. The processes $\wh X$
and $\wh Y$ are defined on the same probability space but no
further relationship such as independence is assumed. In
particular, the two processes do not necessarily reach $z$ at the
same time.

We will now define a distribution for a pair of processes $(\wt
X,\wt Y)$ starting from $x,y \in \gamma_{n-m_0}\cup\gamma_n$, such
that either $x=y$ or $x \in \gamma_{n-m_0}$ and $y \in \gamma_n$.
If $x=y \in \gamma_{n-m_0}$, we define  $\wt X_t = \wt Y_t $ for
all $t\geq 0$, and the distribution of $\wt X$ is that of
reflecting Brownian motion in $D$, killed upon hitting $\gamma_n$.
Similarly, if  $x=y \in \gamma_n$, define $\wt X_t = \wt Y_t $ for
all $t\geq 0$, and the distribution of $\wt X$ is that of
reflecting Brownian motion in $D$, killed upon hitting
$\gamma_{n-m_0}$, The most significant case is when $x \in
\gamma_{n-m_0}$ and $y \in \gamma_n$. In this case we let
$$Z_x \df \eta_{x,n+1} I_{x,n+1} + \eta_{y, n-m_0-1} I_{y,
n-m_0-1} + \eta_{x,y,n-m_0-1} I_{x,y}
$$
and
$$ Z_y \df \eta_{x,n+1} I_{x,n+1} +
\eta_{y, n-m_0-1} I_{y, n-m_0-1} + \eta_{y,x,n+1} I_{x,y}.
$$
Note that by our construction, $Z_x$ and $Z_y$ have
distributions
 $\mu_x$ and $\mu_y$, respectively. When $ I_{y, n-m_0-1}=1$ and
$Z_x=Z_y= \eta_{y, n-m_0-1} =z\in \gamma_{n-m_0-1}$, we define the
(conditional) distribution of $(\wt X,\wt Y)$
to
 be $\Q_{x,y}^z$ until the processes hit $\gamma_{n-m_0-1}$, and
then we ``continue them as a single reflecting Brownian motion in
$D$ starting from $z$ until it hits $\gamma_{n-m_0}$.'' In other
words, if $\wt X$ hits $\gamma_{n-m_0-1}$ at time $t_0$ and $\wt
Y$ hits $\gamma_{n-m_0-1}$ at time $t_1$ then $\wt X_{t_0+t} = \wt
Y_{t_1+t}$ for $t\geq 0$. Similarly, when $I_{x,n+1}=1$ and
$Z_x=Z_y= \eta_{x,n+1}=z\in \gamma_{n+1}$, we define the
(conditional) distribution of $(\wt X,\wt Y)$
to
 be $\Q_{x,y}^z$ until the processes hit $\gamma_{n+1}$, and then
we continue them as a single reflecting Brownian motion in $D$
starting from $z$ until it hits $\gamma_n$. When $I_{x,y}=1$ and
$Z_x =z_1\in \gamma_{n-m_0-1}$ and $Z_y=z_2\in \gamma_{n+1}$, we
let  $\wt X$ have (conditional) distribution $\Q^x_{z_1}$ and
then we continue it as reflecting Brownian motion in $D$
starting from $z_1$ until it hits $\gamma_{n-m_0}$, and we let
$\wt Y$ be independent from $\wt X$ with conditional
distribution $\Q_y^{z_2}$, and we continue it as reflecting
Brownian motion in $D$ starting from $z_2$ until it hits
$\gamma_n$. We call the distribution of the processes constructed
above $\P_{x,y}$. Note that under $\P_{x,y}$ each one of the
processes $\wt X$ and $\wt Y$ is a reflecting Brownian motion in
$D$. Under
$\P_{x,y}$,
 the processes $\wt X$ and $\wt Y$ start from $x,y\in
\gamma_{n-m_0}\cup \gamma_n$, i.e., $\wt X_0=x$ and $\wt Y_0=y$,
they have random lifetimes $\zeta^X$ and $\zeta^Y$, not
necessarily equal, $\wt X_{\zeta^X-}, \,  \wt Y_{\zeta^Y-} \in
\gamma_{n-m_0}\cup \gamma_n$, and either $\wt X_{\zeta^X-}= \wt
Y_{\zeta^Y-}$ or $\wt X_{\zeta^X-} \in \gamma_{n-m_0}$ and $\wt
Y_{\zeta^Y-} \in \gamma_n$. The essential property of $\P_{x,y}$
is that if $\wt X$ enters $\Omega_{n+1}$ before it is killed, then
the part of the trajectory of $\wt X$ after the hitting time of
$\Omega_{n+1}$ is a time shift of the trajectory of $\wt Y$ after
its hitting time of $\Omega_{n+1}$. Similarly, under $\P_{x, y}$,
if $\wt Y$ enters $\Omega'_{n-m_0-1}$ before it gets killed, then
the part of the trajectory of $\wt Y$ after the hitting time of
$\Omega_{n-m_0-1}$ is a time shift of the trajectory of $\wt X$
after its hitting time of $\Omega_{n-m_0-1}$.

We will use the distributions $\P_{x,y}$ to construct processes
$X$ and $Y$ which are defined on  the whole time interval $[0,\infty)$.
Suppose that $x \in
\gamma_{n-m_0}$ and $y \in \gamma_n$ and let $(X^1,Y^1)$ have
distribution $\P_{x,y}$. Let $(X^2,Y^2)$ have conditional
distribution $\P_{x_2,y_2}$ given the event $\big\{X^1_{\zeta^{X^1}-}=x_2
\hbox{ and } Y^1_{\zeta^{Y^1}-}=y_2\big\}$. We continue by
induction. Given $(X^k,Y^k)$, we let $(X^{k+1},Y^{k+1})$ have
conditional distribution $\P_{x_{k+1},y_{k+1}}$ given
the event $\big\{X^k_{\zeta^{X^k}-} = x_{k+1} \hbox{ and }
Y^k_{\zeta^{Y^k}-} = y_{k+1}\big\}$. It is easy to see that
$\sum_k \zeta^{X^k} =\infty$ and $\sum_k \zeta^{Y^k} =\infty$,
a.s. Set $\zeta^{X^0} = \zeta^{Y^0} =0$. For $k\geq 0$ and $t\in
\big[\sum_{0\leq j\leq k} \zeta^{X^j}, \sum_{0\leq j\leq k+1}
\zeta^{X^j}\big)$, define
$$X_t \df X^{k+1}\Big( t- \sum_{0\leq j\leq k} \zeta^{X^j}\Big).
$$
Similarly, for $t\in \big[\sum_{0\leq j\leq k} \zeta^{Y^j},
\sum_{0\leq j\leq k+1} \zeta^{Y^j}\big)$, define
$$ Y_t \df Y^{k+1}\Big(t- \sum_{0\leq j\leq k} \zeta^{Y^j}\Big).
$$
It is straightforward to check that $X$ and $Y$ are reflecting
Brownian motions in $D$ and $\{T^X_K < T^X_{B_*}\} \subset
\{T^Y_K<T^Y_{B_*}\}$. This proves that $\P_x(T_K < T_{B_*}) \leq
\P_y(T_K < T_{B_*})$ for $x \in \gamma_{n-m_0}$ and $y \in
\gamma_n$ and, as we pointed out at the beginning of this proof,
this implies that $\sup_{x\in \gamma_{n-m_0}} G(x_0,x) \leq
\inf_{x\in \gamma_n} G(x_0,x)$.
 \qed

Recall $\lambda$ from Definition 2.1.

\bigskip
\noindent{\bf Lemma 4.11}. {\sl
 For any $c_1$ and $\lambda$ there exists $c_2$ such that the
following holds. Suppose that for some $n$ and $m_2$ we have
$|D_k| \leq c_1 r_k^d$ for all $n-m_2 -1\leq k \leq n$. Then for
all $x \in \gamma_n$,
 $$\P_x (T_{\gamma_{n-m_2-1}} < T_{\gamma_{n+1}})
 \leq c_2  {r_{n}^{2-d}\over\sum_{i=n-m_2 -1}^{n} r_i^{2-d}},$$
and for all $x \in \gamma_{n-m_2}$ we have
 $$\P_x(T_{\gamma_{n+1}} < T_{\gamma_{n-m_2-1}})
 \leq c_2  {r_{n-m_2 -1}^{2-d}\over\sum_{i=n-m_2 -1}^{n} r_i^{2-d}}.$$

}

\bigskip
\noindent {\bf Proof}. We prove the second inequality, the first
one being very similar.
Write $j$ for $n-m_2-1$. If $a,b$ are integers
with $j\leq a\leq b\leq n$, set $U_{a,b}=\bigcup_{k=a}^b D_k$,
define
 $$C_{a,b}=\inf\Bigl\{ \int_{U_{a,b}} |\nabla f(x)|^2 dx: f=0\hbox{ on }
\gamma_a\hbox{ and }
f=1 \hbox{ on }\gamma_{b+1}\Bigr\}, \eqno (4.16)$$ and let
$R_{a,b}=C_{a,b}^{-1}$. $C_{a,b}$ is called the conductance across
$U_{a,b}$ and $R_{a,b}$ the resistance. Consider reflecting
Brownian motion in $U_{a,b}$ killed on hitting $\gamma_a$ and let
$G_{a,b}(x,y)$ be the corresponding Green function. We use the
fact that with respect to this process $C_{a,b}$ is equal to the
capacity of $\gamma_{b+1}$; see [FOT].

Using Definition 4.1 we can find a constant $c_3$ independent of
$k$ and points $z_k\in D_k$ such that $\,\hbox{dist}\,(z_k,
\partial D_k)\geq c_3 r_k$. Let $B_k$ be the ball of radius
$c_3r_k/2$ centered at $z_k$. Starting at any point that is a
distance $c_3r_k/4$ from $z_k$, the expected time that Brownian
motion in $D_k$ spends in $B_k$ before hitting $\partial D_k$ is
larger than $c_4 r_k^2$. By the support theorem for standard
$d$-dimensional Brownian motion, starting from any point that is a
distance $3c_3r_k/4$ from $z_k$, there is probability at least
$p_1>0$ (not depending on $k$) that the Brownian motion will hit
the ball of radius $c_3r_k/4$ about $z_k$ before hitting $\partial
D_k$.  So starting at such a point the expected time spent in
$B_k$ before hitting $\partial D_k$ is at least $p_1 c_4 r_k^2$.
Using the Harnack inequality and the fact that $|B_k|=c_5 r_k^d$,
it follows that $G_{k,k}(z_k,y)\geq c_6r_k^{2-d}$ if
$|y-z_k|=3c_3r_k/4$. By the Neumann boundary Harnack principle,
$$G_{k,k}(z_k,y)\geq c_7 r_k^{2-d}, \qquad y\in \gamma_{k+1}. \eqno (4.17)$$

Consider reflecting Brownian motion in $D_k$ killed on hitting
$\gamma_k$ and let $\nu_k$ be the capacitary measure for
$\gamma_{k+1}$. Then
$$\eqalign{1&\geq {\bf P}_{z_k}(T_{\gamma_{k+1}}<T_{\gamma_k})
=\int G_{k,k}(z_k,y) \nu_k(dy)\cr &\geq c_7 r_k^{2-d}
\nu_k(\gamma_{k+1})\cr &= c_7 r_k^{2-d} C_{k,k}.\cr}$$ Therefore
$$C_{k,k}\leq c_7^{-1} r_k^{d-2}$$
and
$$R_{k,k}\geq c_7 r_k^{2-d}. \eqno (4.18)$$

Next, if $a_1\leq a_2<a_2+1\leq a_3$, let $f_1$ be the function on
$U_{a_1,a_2}$ at which the infimum in (4.16) is attained and
similarly $f_2$ the function on $U_{a_2+1,a_3}$.
Let $\beta=C_{a_2+1,a_3}/(C_{a_1,a_2}+C_{a_2+1,a_3})$
 and define $f$ on $U_{a_1,a_3}$ by setting the restriction of $f$
on $U_{a_1,a_2}$ to be equal to $\beta f_1$ and the restriction of
$f$ on $U_{a_2+1,a_3}$ to be equal to $\beta+(1-\beta)f_2$. Then
$$\eqalignno{ C_{a_1,a_3}&\leq \int_{U_{a_1,a_3}}
|\nabla f(x)|^2 dx=\beta^2\int_{U_{a_1,a_2}}|\nabla f_1|^2
+(1-\beta)^2 \int_{U_{a_2+1,a_3}} |\nabla f_2|^2\cr &=\beta^2
C_{a_1,a_2}+(1-\beta)^2 C_{a_2+1, a_3}\cr
&={{C_{a_1,a_2}C_{a_2+1,a_3}}\over {C_{a_1,a_2}+C_{a_2+1,a_3}}}.
&(4.19)\cr}$$ This is equivalent to
$$R_{a_1,a_3}\geq R_{a_1,a_2}+R_{a_2+1,a_3}. \eqno (4.20)$$

By (4.18), (4.20), and induction, we obtain
$$R_{j,n}\geq \sum_{i=j}^n c_7 r_i^{2-d},$$
or
$$C_{j,n}\leq {1\over{\sum_{i=j}^n c_7 r_i^{2-d}}}. \eqno (4.21)$$

Recall that $B_j$ is the ball of radius $c_3r_j/2$ about $z_j$.
Starting in $B_j$ the expected amount of time the process spends
in $B_j$ before exiting the ball of radius $3c_3r_j/4$ about $z_j$
is bounded by $c_8r_j^2$. By the support theorem for standard
Brownian motion, there exists $p_2>0$ such that starting at any
point that is a distance $3c_3r_j/4$ from $z_j$, there is
probability at least $p_2$ of hitting $\gamma_j$ before hitting
$\partial D_j\setminus \gamma_j$. A standard argument allows us to
conclude that the expected amount of time spent in $B_j$ starting
at any point of $U_{j,n}$ is at most $c_9 r_j^2$. Since
$|B_j|=c_{10}r_j^d$, the Harnack inequality implies that
$G_{j,n}(z_j,y)\leq c_{11} r_j^{2-d}$ if $y\in \gamma_{n+1}$. The
Neumann boundary Harnack inequality then implies that
 $$G_{j,n}(x,y)\leq c_{12}r_j^{2-d}, \qquad x\in \gamma_{j+1},
 \quad y\in \gamma_{n+1}. \eqno (4.22)$$
Let $\nu$ be the equilibrium measure for $\gamma_{n+1}$ with
respect to reflecting Brownian motion in $U_{j,n}$ killed on
hitting $\gamma_j$. Combining (4.21) and (4.22),
$${\bf P}_x(T_{\gamma_{n+1}}<T_{\gamma_j})=\int G_{j,n}(x,y) \nu(dy)
\leq c_{12} r_j^{2-d} C_{j,n}, \qquad x\in \gamma_{j+1}.
$$
This proves the lemma. \qed

\bigskip

If the $r_n$ are comparable, then Lemma 4.11 implies Condition 4.7
 for sufficiently large $m_0$.

\bigskip
\noindent{\bf Remark 4.12}. If $D\in \cD_1$, then a ``typical
point'' $x\in \prt D$ has a neighborhood $U\subset \ol D$ such
that $\prt D\cap U$ is the graph of a Lipschitz function. The
Green function satisfies $G(x,y) \leq c_1 |x-y|^{2-d}$, for $x,y
\in U$, where $c_1$ depends only on the Lipschitz constant
characterizing $\prt D$; to see this we flatten the boundary and
reflect over a hyperplane as in the proof of Lemma 2.8, and then
use the result of [LSW]. The upper estimate in Lemma 4.4 follows
from this immediately.

\bigskip

\noindent{\bf Example 4.13}. Our first example in this section is a
multidimensional version of Example 3.4. Suppose that $d\geq 3$
and for some $\alpha>1$,
 $$
 D= \left\{x = (x_1, x_2, \dots, x_d): 0 < x_1 < 1
\hbox{ and } x_1^\alpha >(x_2^2 + \dots + x_d^2)^{1/2}\right\}.
 $$
We will restrict the parameter range to $\alpha>1$. We will show
that if $\alpha \in (1,2)$ then the whole surface of $D$ is active
and when $\alpha \geq 2$ then part of the surface is nearly inactive.

We will analyze only one boundary point, the origin, in view of
Remark 4.12. We let the $\gamma_k$'s be intersections of $D$ with
 $(d-1)$-dimensional hyperplanes perpendicular to the first axis,
at distances $2^{-k} + j 2^{-k\alpha}$ from 0, for all $j\geq 0$
such that $2^{-k} + j 2^{-k\alpha}\leq 2^{-k+1} - 2^{-k\alpha}$,
for all $k\geq 1$.

Note that for some $c_1$ and any $m_0$ there exists $m_1$ such
that for any $n>m_1$ we have $1/c_1 \leq r_j/r_k \leq c_1$ for all
$n\leq j,k\leq n+m_0$. This and Lemma 4.11 easily imply that
Condition  4.7 holds.

The number of $D_n$'s whose distance from 0 lies between $2^{-k}$
and $2^{-k+1}$ is of order $2^{-k(1-\alpha)}$. For $D_n$'s in this
range, $\sum_{m=1}^n r_m^{2-d}\approx \sum_{j\leq k}
2^{-j(1-\alpha)} 2^{-j\alpha(2-d)}\approx 2^{-k(1+\alpha(1-d))}$.
The surface area, $|\prt D_n \cap \prt D|$, is of order $2^{-k
\alpha(d-1)}$, so the contribution from these sets to the sum in
(4.1) is of order $2^{-k(1-\alpha)}\cdot 2^{-k(1+\alpha(1-d))}
\cdot 2^{-k \alpha(d-1)} = 2^{-k (2-\alpha)}$. If $\alpha < 2$
then $\sum_{k\geq 1} 2^{-k (2-\alpha)}< \infty$, so part (i) of
Theorem 4.3 implies that the whole surface of $D$ is active.

A similar calculation shows that the sum in (4.2) is comparable to
$\sum_{k\geq 1} 2^{-k (2-\alpha)}$ and this is infinite for
$\alpha \geq 2$. Hence, by Theorem 4.3 (ii), part of the surface
of $D$ is nearly inactive if $\alpha \geq 2$.  \qed

\bigskip
\noindent{\bf Remark 4.14}. Fukushima and Tomisaki [FT] studied
reflecting Brownian motion in unbounded cusps
$$ \wt D\df \left\{x = (x_1, x_2, \dots, x_d): x_1>0 \hbox{ and }
 x_1^\alpha > (x_2^2 + \dots + x_d^2)^{1/2}  \right\}
 $$
and derived a Green function estimate (see Lemma 5.4 and 5.5 in
[FT]). Their proof can be adapted to get the Green function upper
bound estimate for reflecting Brownian motion in the truncated
cusps $D$ as defined in Example 4.13  and to show that for
$1<\alpha<2$,
$$ \int_{\partial D} G_{D\setminus B_*} (x, y)
\sigma (dy)<\infty.
$$
Thus this gives an alternative proof for the boundary of $D$ to be
active when $1<\alpha<2$. The main goal of the paper [FT] is to
show that reflecting Brownian motion in $\wt D$ starting from the
cusp point ${\bf 0}\df (0, \cdots, 0)$ is a semimartingale when
$\alpha <2$. Using Theorem 4.3(ii) (and its proof) in this paper,
we can settle the remaining case by showing that reflecting
Brownian motion in $\wt D$ starting from ${\bf 0}$ is not a
semimartingale when $\alpha \geq 2$. Clearly, this is equivalent
to the fact that reflecting Brownian motion in $D$ starting from
${\bf 0}$ is not a semimartingale when $\alpha \geq 2$.

By Theorem 2.1 of [FT], reflecting Brownian motion $X$ in $D$ is a
strong Feller process on $\overline D$ and thus can start from
every point in $\overline D$. Let $\alpha \geq 2$. According to
Example 4.13, part of $\partial D$ is nearly inactive. By the
proof of Theorem 4.3(ii),
$$
\lim_{x\to {\bf 0}} \P_x \left( L_{T_{B_*}}>b \right)>1-\frac1{b}
\qquad \hbox{for every } b>0.
$$
Were $X$ a semimartingale starting from ${\bf 0}$, the Skorokhod
decomposition for $X$
$$ X_t=X_0+W_t+\int_0^t \n (X_s) dL_s \qquad \hbox{for } t\geq 0
$$
would hold under $\P_x$ for every $x\in \overline D$. It follows
from weak convergence and the second to the last display that
$$ \P_{\bf 0} (L_{T_{B_*}} =\infty )=1.
$$
This is a contradiction since $\P_{\bf 0} (T_{B_*} <\infty) >0$.
Therefore $X$ starting from ${\bf 0}$ cannot be a semimartingale.
\qed

\bigskip
\noindent{\bf Example 4.15}. This is a multidimensional analogue
of Example 3.6. Suppose that $\alpha>0$, $\beta >1$ and let $a_k =
\sum_{j=1}^k 2^{-(j-1)\alpha}$. Let $\cS_n$ be the family of all
binary (zero-one) sequences of length $n$. We will write $\s =
(s_1, s_2, \dots, s_n)$ for $\s\in \cS_n$. For integer $k\geq 1$
and $\s\in \cS_k$, we set $b_{\s} = \sum_{j=1}^k s_j 2^{-j}$. Let
$A_* = [0,1]^d$. For $k\geq 1$ and $\s\in \cS_k$ let
 $$A_{\s} = \{(x_1,\dots,x_d): a_k \leq x_1\leq a_{k+1},
((x_2-b_\s)^2+x_3^2+\dots+x_d^2)^{1/2} \leq f_\s(x_1)
 \},$$
where $c_1 2^{-k\beta} \leq f_\s(x_1) \leq 2^{-k\beta}$ and
$c_1>0$ does not depend on $\s$. Assume that all functions $f_\s$
are Lipschitz with the same Lipschitz constant. Let $D$ be the
connected component of the interior of $A_*\cup \bigcup_{k\geq 1}
\bigcup_{\s\in\cS_k} A_{\s}$ that contains the open box $(0,
1)^d$.

We restrict the range of parameters to $\beta >\alpha$. Since we
have assumed that $\alpha>0$ and $\beta>1$, the surface of $D$ is
finite and Remark 3.5 cannot be used to draw any conclusions.

We will show that if $\alpha < \beta < 2\alpha$ then the whole
surface of $D$ is active and when $\beta \geq 2\alpha$ then
part of the surface is nearly inactive.

We will analyze only a family of $D_n$'s corresponding to a
boundary point at the end of a channel, in view of Remark 4.12.
 Fix a boundary point $z_0$ at the end of an infinite channel,
i.e., a point whose first coordinate is $\sum_{j=1}^\infty
2^{-(j-1)\alpha}$. Let ${\cal A}_k$ be the family of hyperplanes
$K_{k,n}=\{(x_1,\dots,x_d): x_1=a_k + n 2^{-k \beta}\}$, with
$n\geq 1$ such that $a_k + n 2^{-k \beta} \leq a_{k+1} $. Let
${\cal C}_k$ be the family of connected components of $K_{k,n}
\cap D$, for $K_{k,n} \in {\cal A}_k$, which separate $z_0$ from
$A_*$. Let $\gamma_n$'s be the relabelled family $\bigcup_k {\cal
C}_k$.

This assumption on the magnitude of $f$ and Lemma 4.11 imply that
Condition  4.7 holds as long as the relevant $D_n$'s belong
to the same $A_{\s}$.
 Lemmas 4.9 and 4.10 then prove that Condition 4.2, i.e.,
$\sup_{x\in \gamma_{n-m_0}} G(x_0,x) \leq \inf _{x\in \gamma_{n}}
G(x_0,x)$, holds if
$\gamma_{n-m_0}$ and $\gamma_{n}$
 belong to the same
$A_{\s}$.
 Then clearly Condition  4.2 holds in full generality if we
replace $m_0$ with $2m_0 +1$.

The number of $D_n$'s defined by the $\gamma_n$'s, needed to reach
$A_\s$ with $\s\in\cS_k$ is of order $\sum_{j\leq k}
2^{-j\alpha}/2^{-j \beta}\approx 2^{k(\beta-\alpha)}$. Consider a
$D_n$ which intersects $A_\s$ with $\s\in\cS_k$. The set $D_n$ may
either have diameter of order $2^{-k \beta}$ or it may contain
a ``tree'' of thin channels. Consider first $D_n$'s that have
diameters of order $2^{-k \beta}$. There are $2^{k(\beta-\alpha)}$
such $D_n$'s, up to a constant, so for $n$ in this range,
$\sum_{m=1}^n r_m^{2-d}\approx \sum_{j\leq k} 2^{j(\beta-\alpha)}
2^{-j\beta(2-d)}\approx 2^{k(\beta(d-1)-\alpha)}$. The surface
area, $|\prt D_n \cap \prt D|$, is of order $2^{-k\beta(d-1)}$, so
the total contribution of such $D_n$'s to (4.1) is of order
$2^{k(\beta-\alpha)}2^{k(\beta(d-1)-\alpha)}2^{-k\beta(d-1)}
\approx 2^{k(\beta -2\alpha)}$. The series $\sum_k 2^{k(\beta
-2\alpha)}$ is summable if and only if $\beta < 2 \alpha$.

Next consider a $D_n$ which intersects $A_\s$ with $\s\in\cS_k$
and contains a side ``tree'' of thin channels. Its surface area,
$|\prt D_n \cap \prt D|$, is of order $\sum_{j\geq k} 2^{j-k}
2^{-j\alpha} 2^{-j\beta(d-2)}\approx 2^{-k(\alpha+\beta(d-2))}$.
There are at most two such $D_n$'s for each $A_\s$, so their
contribution to (4.1) is of order $2^{k(\beta(d-1)-\alpha)}
2^{-k(\alpha+\beta(d-2))}\approx 2^{k(\beta -2\alpha)}$. Hence,
the contribution of $D_n$'s with side channels is of the same
order as the contribution of $D_n$ that have diameter of order
$2^{-k \beta}$. We conclude that (4.1) holds if $\beta < 2
\alpha$.

If $\beta \geq 2 \alpha$ then the contribution of $D_n$'s with
diameter of order $2^{-k \beta}$ is enough to make the left hand
side of (4.2) infinite, due to the estimates presented above.
\qed
\bigskip

{\bf Example 4.16}. We will analyze a multidimensional fractal
domain vaguely resembling the von Koch snowflake, except that we
will add barriers partly blocking the passage between the building
blocks. Suppose that the dimension of the space is $d\geq 3$ and
fix a parameter $\rho \in (0,1/2)$.
We will impose further restrictions on $\rho$ below.
 For $k\geq 0$, let $\cA_k$ be a finite family of open cubes with
edge length $\rho^k$, with edges parallel to the axes, and
satisfying the following properties. The family $\cA_0$ consists
of one cube $A_0$. The family $\cA_1$ consists of $2d$ cubes which
are disjoint from each other and are disjoint from $A_0$. One side
of any cube in $\cA_1$ lies on a side of $A_0$ and these two sides
of the two cubes have the same center. Now suppose that we have
defined families $\cA_k$ for $k\leq n$. Let $A_n$ be the union of
all cubes in $\bigcup _{k\leq n} \cA_k$. Then $\cA_{n+1}$ is the
maximal family of disjoint cubes that do not intersect $A_n$ and
such that one side of each of these cubes lies on a side of a cube
from the family $\cA_n$, and has the same center. Let $D_*$ be the
union of all cubes in $\bigcup _{k\geq 0} \cA_k$ and note that
this set is not connected because all cubes in this family are
disjoint. We transform $D_*$ into a connected open set by adding
``passages'' between cubes. Fix a parameter $\beta>1$. For any
pair of cubes which belong to $\cA_{n-1}$ and $\cA_{n}$, and whose
sides intersect and have a common center $x$, we add to $D_*$ the
open ball $B(x,2\rho^{(n-1) \beta})$. We let $D$ be the union of
$D_*$ and all such balls. Parts of the boundary of $D$ are
adjacent to $D$ on both sides, and this is forbidden by Definition
4.1, strictly speaking. We could modify the domain $D$ or even
Definition 4.1 to cover this case, but that would be an
unnecessary embellishment.

We will determine for which values of $\rho$ the surface area is
finite because the example is not interesting if $|\prt
D|=\infty$; in such a case a part of the surface is nearly
inactive by Remark 3.5. The surface area of a cube with edge
length $\rho^k$ is of order $\rho^{k(d-1)}$. The number of cubes
in $\cA_k$ is of order $(2d-1)^k$. The total surface area of cubes
in $\cA_k$ is of order $\rho^{k(d-1)}(2d-1)^k$. The surface area
of $D$ is finite if $\sum_k \rho^{k(d-1)}(2d-1)^k < \infty$, that
is if $\rho^{d-1} (2d-1) < 1$. Hence, we are interested only in
$\rho$ less than $ (2d-1)^{-1/(d-1)}$. The function $
f(d)=(2d-1)^{-1/(d-1)}$ is increasing for $d\geq 3$ because, when
we treat $d$ as a real argument,
 $$ f'(d) =\frac{1 + \left(  2d -1\right) ( \log (
 2d-1)-1)}
  {{\left( d -1\right) }^2\,{\left( 2d -1\right) }^{d/( d-1)}}
  >0
 $$
for $d\geq 3$. We have $ f(3) = 1/\sqrt{5}$, and $\lim_{d\to
\infty} f(d) = 1$. Hence, we can take $\rho\in(0,1/2\land
1/\sqrt{5})$ for any $d$. We will see that, as long as the surface
area is finite, the value of $\rho$ does not play any role in this
example.

As usual in our examples, we will analyze only a point $z\in\prt
D$ that lies at the end of an ``infinite'' channel, i.e., such
that any continuous path in $D$ from the center $z_*$ of $A_0$ to
$z$ must pass through at least one cube in every family
$\cA_k$. Let $\Gamma$ be a continuous path in $D$ from $z_*$ to
$z$ that passes through a side of any cube in $\bigcup
_{k\geq 0} \cA_k$ at most once, and if it does so, then it passes
through the center of that side. Let $z_k$ be the intersection
point of $\Gamma$ and the side of the cube in $\cA_k$ that is
a part of a side of a cube in $\cA_{k-1}$. The curve $\Gamma$
passes through all the $z_k$'s on its way from $z_*$ to $z$.

For every $z_k$, let $\cC_k$ be the family of all sets $\prt
B(z_k, 2^j  ) \cap D$, where $j$ satisfies $4\rho^{(k-1)\beta}
\leq 2^{j-1} \leq 2^{j+1} \leq \rho^{k-1} /2$. Note that each set
$\prt B(z_k, 2^j  ) \cap D$ contributes two sets to $\cC_k$ and
each one of these sets is a spherical cap. Let $\cC =\bigcup_k
\cC_k$ and rename the elements of $\cC$ as $\gamma_n$, in the order in
which they have to be passed on the way from $z_*$ to $z$ within
$D$. It is elementary to check that this family of $\gamma_n$'s
satisfies the conditions listed in Definition 4.1.

In this example, Condition  4.7 does not hold. We will argue that
Condition 4.8 holds directly. Consider spheres $\prt B(0, 2^{j}), \prt
B(0, 2^{j+k_0}), \prt B(0, 2^{j+k_0+m_0})$ and $\prt B(0,
2^{j+2k_0+m_0})$ and call them $S_1, S_2, S_3$ and $S_4$. It is
not very hard to prove that there exist large $k_0$ and  $m_0$,
such that $\P_x(T_{S_4} \in A, T_{S_4} < T_{S_1}) \leq
\P_y(T_{S_4} \in A, T_{S_4} < T_{S_1})$, for $A\subset S_4$, $x\in
S_2$ and $y\in S_3$. We also have, for sufficiently large $k_0$
and $m_0$, that $\P_x(T_{S_1} \in A, T_{S_1} < T_{S_4}) \geq
\P_y(T_{S_1} \in A, T_{S_1} < T_{S_4})$, for $A\subset S_1$, $x\in
S_2$ and $y\in S_3$, although the two claims are not symmetric and
require somewhat different justification. By the reflection
principle, for reflecting Brownian motion in $D$,
$\P_x(T_{\gamma_{n+2k_0+m_0}} \in A, T_{\gamma_{n+2k_0+m_0}} <
T_{\gamma_n}) \geq \P_y(T_{\gamma_{n+2k_0+m_0}} \in A,
T_{\gamma_{n+2k_0+m_0}} < T_{\gamma_n})$ for $A\subset
\gamma_{n+2k_0+m_0}$, $x\in \gamma_{n+k_0}$ and $y \in
\gamma_{n+k_0+m_0}$, provided $\gamma_{n}$ and
$\gamma_{n+2k_0+m_0}$ belong to the same family $\cC_k$ and lie on
the same side of $\prt D$. We also have $\P_x(T_{\gamma_{n}} \in
A, T_{\gamma_{n}} < T_{\gamma_{n+2k_0+m_0}}) \geq
\P_y(T_{\gamma_{n}} \in A, T_{\gamma_{n}} <
T_{\gamma_{n+2k_0+m_0}})$ for $A\subset \gamma_{n}$, $x\in
\gamma_{n+k_0}$ and $y \in \gamma_{n+k_0+m_0}$. If we take only
every $k_0$-th element of the family $\gamma_n$, this proves
Condition 4.7 for $\gamma_n$'s which belong to the same family $\cC_k$ and
lie on the same side of $\prt D$. Hence, Lemma 4.10 proves
Condition  4.2 for $n$ restricted in such a way. However,
this implies that Condition  4.2 holds for all $n$, with
$m_0$ replaced by $2m_0$, for the same reason as in Example 4.13.

The number of $\gamma_n$'s in $\cC_k$ is of order $\log
((\rho^{k-1} /2)/( 4\rho^{(k-1)\beta})) \approx k$. If $z_1\in
\Omega_{n+1}$ with sufficiently large $n$ then by Lemma 4.4, the
Green function $G(z_1, \,\cdot\,)$ can be bounded by $c_1
\sum_{j\leq k} j \rho^{j \beta (2-d)} \leq c_2 k \rho^{k \beta
(2-d)}$ for $x\in D$ that lie between $\gamma_n$'s in $\cC_k$. The
surface area of $\prt D_n\cap \prt D$ corresponding to $\gamma_n
\in \cC_k$ is bounded by $c_3 \rho^{k(d-1)}$, so the contribution
of such $D_n$'s to the sum in (4.1) is bounded by $c_4 k^2
\rho^{k\beta(2-d)} \rho^{k(d-1)} = c_4 k^2 \rho^{k(\beta(2-d)
+d-1)}$. If $\beta < (d-1)/(d-2)$ then $\sum_k
k^2\rho^{k(\beta(2-d) +d-1)}<\infty$ and Theorem 4.3 (i) implies
that the whole surface of $D$ is active.

To find a lower bound for (4.2), we take into account only one
$D_n$ corresponding to each family $\cC_k$, namely the one with
the largest surface area. We obtain as a  lower bound for
(4.2) the quantity  $c_5 k\rho^{k\beta(2-d)} \rho^{k(d-1)} = c_5
k\rho^{k(\beta(2-d) +d-1)}$. If $\beta \geq (d-1)/(d-2)$ then
$\sum_k k\rho^{k(\beta(2-d) +d)}=\infty$ so by Theorem 4.3 (ii),
part of the  surface of $D$ is nearly inactive.

\bigskip

\noindent{\bf 5. Trap domains.}

The ideas developed in Section 4 allow us to prove a new result on
``trap'' domains introduced in [BCM]. The new result applies only
to the class of domains $\cD_1$ presented in Definition 4.1 but
that class contains some very natural examples of fractal domains,
such as the multidimensional version of the von Koch snowflake
presented in Example 4.16, that were not covered by theorems
proved in [BCM] (see Example 5.2 below). There was a big gap
between results on two-dimensional domains and higher dimensional
domains in [BCM]. At first we thought the gap was purely technical
in nature---complex analytic methods could not be used in higher
dimensions. It turns out that the gap is in fact ``real,'' in the
sense that the multidimensional examples are considerably
different from the two-dimensional examples---compare our Example
5.2 and Proposition 2.15 of [BCM].

Recall that $B_*\subset D$ is a closed ball with positive radius
and $T_{B_*} = \inf\{t\geq 0: X_t \in B_*\}$ is the first hitting
time of $B_*$ by $X$. We say that $D\subset \R^d$, $d\geq 2$, is a
{\it trap domain} if
$$\sup_{x\in D} \E_x T_B = \infty,\eqno(5.1)
$$
and otherwise $D$ is called a {\it non-trap domain}. One can
express (5.1) in a purely analytic way, namely, by saying that $D$
is a trap domain if and only if
$$\sup_{x\in D\setminus B} \int_{D\setminus B} G(x,y) dy
= \infty. \eqno(5.2)
$$

See [BCM] for further discussion of basic properties of trap
domains.

\bigskip

\noindent{\bf Theorem 5.1}. {\sl Consider
 a domain $D\in \cD_1$, $D\subset \R^d$, $d\geq 3$, with a finite
volume.

(i) If there exists a
 constant $c<\infty$ such that for each
 point $z\in \b D$, there is a system of surfaces $\{\gamma_n,
n\geq 0\}$ as in Definition 4.1 satisfying
$$  \sum_{n=1}^\infty |D_n| \sum_{k=1}^n r_k^{2-d}
  \leq c,\eqno(5.3)
$$
then $D$ is not a trap domain.

(ii) If there exists a boundary point $z\in \b D$ and a system of
surfaces $\{\gamma_n\}$ as in Definition 4.1, such that
 $$ \sum_n |D_n| \sum_{k=1}^n r_k^{2-d}=\infty,\eqno(5.4)$$
then $D$ is a trap domain.

}

\bigskip

\noindent{\bf Proof}. (i) The proof is very similar to the proof
of Theorem 4.3. Consider $x_0 \in D$. It is not hard to see that
there exists $z_0\in \prt D$ and a corresponding family of
$\gamma_n$'s, as in Definition 4.1, such that $x_0 \in D_{n_0}$
for some $n_0$ and $\dist(x_0, \prt D_{n_0}) \geq c_1 r_{n_0}$. By
Lemma 4.4, $G(x_0, \,\cdot \,)$ is bounded by $c_2
\sum_{k=1}^{n_0} r_k^{2-d}$ on $D_{n_0-1}$. By the Harnack
principle, it is bounded by $c_3 \sum_{k=1}^{n_0} r_k^{2-d}$ on
$\prt B(x_0, c_1 r_{n_0}/2)$, and the maximum principle implies
that the same bound holds on $D \setminus B(x_0, c_1 r_{n_0}/2)$.
For $x\in B(x_0, c_1 r_{n_0}/2)$ we have $G(x_0,x) \leq c_4
|x-x_0|^{2-d}  \sum_{k=1}^{n_0} r_k^{2-d} / ( r_{n_0}/2)^{2-d}$,
by comparison with the Green function in $\R^d$. We obtain, using
the upper bound in Lemma 4.4 for $n<n_0$,
 $$\eqalign{
 \int_{D\setminus B}& G(x_0,y) dy
 = \sum_{n\geq 1} \int_{D_n} G(x_0,y) dy\cr
 &= \sum_{1\leq n< n_0}\int_{D_n} G(x_0,y)dy
 +\sum_{n\geq n_0}\int_{D_n\setminus B(x_0, c_1 r_{n_0}/2)}
 G(x_0,y)dy\cr
 &\qquad + \int_{ B(x_0, c_1 r_{n_0}/2)} G(x_0,y)dy
 \cr
 &\leq  \sum_{1\leq n< n_0} | D_n |
 c_2 \sum_{k=1}^{n} r_k^{2-d}
 + \sum_{n\geq n_0} | D_n |
 c_3 \sum_{k=1}^{n_0} r_k^{2-d}
 + c_5 r_{n_0}^d  \sum_{k=1}^{n_0} r_k^{2-d}\cr
 &\leq  \sum_{n\geq 1} c_6 | D_n |
 \sum_{k=1}^{n} r_k^{2-d}.
 }$$
This is bounded by a constant independent of $x_0$, by assumption
(5.3). The theorem follows in view of (5.2).

(ii) A calculation similar to that in part (i), based on the lower
bound in Lemma 4.9, easily implies part (ii) of the theorem.
 \qed

\bigskip
We would like to emphasize that part (ii) of Theorem 5.1 is much
easier to prove than part (ii) of Theorem 4.3. This is because all
we have
to show is
 that the function $x\mapsto \E_x \left[T_{B_*}\right]$ is
unbounded. In the proof of Theorem 4.3 (ii) we had to prove that
the random variable $L_{T_{B_*}}$ for reflecting Brownian motion
$X^*$ starting from $x$ converges to infinity in
distribution as $x$ approaches a boundary point $z_0\in
\partial D$.

\bigskip

{\bf Example 5.2}. Recall the domain $D$ and notation from Example
4.15. Recall that the number of $\gamma_n$'s in $\cC_k$ is of
order $\log ((\rho^{k-1} /2)/( 4\rho^{(k-1)\beta})) \approx k$. If
$z_1\in \Omega_{n+1}$ with sufficiently large $n$ then by Lemma
4.4, the Green function $G(z_1, \,\cdot\,)$ can be bounded by $c_1
\sum_{j\leq k} j \rho^{j \beta (2-d)} \leq c_2 k \rho^{k \beta
(2-d)}$ for $x\in D$ that lie between $\gamma_n$'s in $\cC_k$. The
volume of a $D_n$ corresponding to a $\gamma_n \in \cC_k$ is bounded
by $c_3 \rho^{kd}$, so the contribution of such $D_n$'s to the sum
in (5.3) is bounded by $c_4 k^2 \rho^{k\beta(2-d)} \rho^{kd} = c_4
k^2 \rho^{k(\beta(2-d) +d)}$. If $\beta < d/(d-2)$ then $\sum_k
k^2\rho^{k(\beta(2-d) +d)}<\infty$ and Theorem 5.1 (i) implies
that $D$ is not a trap domain.

To find a lower bound for (5.4), we take into account only one
$D_n$ corresponding to each family $\cC_k$, namely the one with
the largest volume. We obtain as a lower bound for (5.4) the quantity
$c_5 k\rho^{k\beta(2-d)} \rho^{kd} = c_5 k\rho^{k(\beta(2-d)
+d)}$. If $\beta \geq d/(d-2)$ then $\sum_k k\rho^{k(\beta(2-d)
+d)}=\infty$ so by Theorem 5.1 (ii), $D$ is a trap domain.
\qed

\bigskip

We will use the above example to compare Theorem 5.1 to a result
about multidimensional trap domains proved in [BCM]. The result in
[BCM] was based on the notion of a $J_\alpha$-domain, used by
Maz'ja in his book on Sobolev spaces [Maz]. Here is an informal
definition of a $J_\alpha$-domain (see [Maz] or [BCM] for the
rigorous definition). We say that $D$ is a $J_\alpha$ domain if
for every smooth $(d-1)$-dimensional surface $\Lambda$ which
divides $D$ into two connected components $D_1^\Lambda$ and
$D_2^\Lambda$, we have $\min(|D_1^\Lambda|, |D_2^\Lambda|) ^\alpha
\leq c_1 |\Lambda|$, where $c_1< \infty$ depends only on $D$ (here
$|\Lambda|$ is the $(d-1)$-dimensional surface area). Theorem 2.4
of [BCM] implies that if $D$ is a $J_\alpha$ domain with $\alpha <
1$ then $D$ is not a trap domain, and there exists a trap domain
$D\in J_1$.

Roughly speaking, one can determine whether the domain $D $ of
Example 5.2 belongs to $J_\alpha$ with a given $\alpha$ by
comparing the surface area of the opening between cubes in
$\cA_k$ and $\cA_{k-1}$ to the volume of the cubes in
$\bigcup_{j\geq k} \cA_j$. The surface area is of order
$\rho^{k\beta (d-1)}$ and the volume is of order $\rho^{k d}$, so
$D$ is a $J_\alpha$ domain with $\alpha < 1$ if $\beta < d/(d-1)$.
Hence, for the family of domains in Example 5.2, Theorem 2.4 of
[BCM] shows that $D$ is not a trap domain if $\beta < d/(d-1)$,
while Theorem 5.1
of this paper
 shows that this holds for all $\beta < d/(d-2)$, and in addition
it shows that this result is sharp. The gap between the power of
the two approaches is not as striking in dimensions $d\geq 3$ as
it is in the 2-dimensional case, discussed in Proposition 2.15 of
[BCM].

\bigskip

\centerline{REFERENCES}
\bigskip

\item{[BBC]} R.~Bass, K.~Burdzy and Z.-Q.~Chen (2005), Uniqueness
for reflecting Brownian motion in lip domains {\it Ann. I. H.
Poincar\'e \bf 41}, 197--235.

\item{[BH]} R. Bass and P. Hsu (1991), Some potential theory for
reflecting Brownian motion in H\"older and Lipschitz domains. {\it
Ann. Probab. \bf 19}, 486--508.

\item{[Bl]} R.M.~Blumenthal (1992), {\it Excursions of Markov
Processes}, Birkh\"auser, Boston.

\item {[Bu]} K.~Burdzy (1987), \it Multidimensional Brownian
Excursions and Potential Theory, \hfil\break
\rm Longman, London.

\item{[BCM]} K.~Burdzy, Z.-Q.~Chen and D.~Marshall (2005), Traps
for reflected Brownian motion. (preprint)

\item{[BT]} K. Burdzy and E. H. Toby (1995), A Skorohod-type lemma
and a decomposition of reflected Brownian motion. {\it Ann.
Probab. \bf 23}, 586--604.

\item{[C1]} Z.-Q. Chen (1993), On reflecting diffusion processes
and Skorokhod decompositions. {\it Probab. Theory Rel. Fields,
\bf 94}, 281-316.

\item{[C2]} Z.-Q. Chen (1996),  Reflecting Brownian motions
and a deletion result for Sobolev spaces of order (1,2).
{\it Potential Analysis, \bf  5}. 383-401.

\item{[CFW]} Z.-Q.~Chen, P.~Fitzsimmons and R.J.~Williams (1993),
Reflecting Brownian motions: quasimartingales and strong
Caccioppoli sets. {\it Potential Anal. \bf 2}, 219--243.

\item{[D]} J.L.~Doob (1984), \it Classical Potential Theory and Its
Probabilistic Counterpart, \break \rm Springer, New York.

\item{[DT]} D.~DeBlassie and E.~Toby (1993), On the semimartingale
representation of reflecting Brownian motion in a cusp. {\it
Probab. Theory Rel. Fields \bf 94}, 505--524.

\item{[FSF]} M.~Felici, B.~Sapoval, and M.~Filoche (2003)
Renormalized random walk study of oxygen absorption in the human
lung, {\it Phys. Rev. Lett. \bf 92}, 068101-1---068101-4.

\item{[FS]} M.~Filoche and B.~Sapoval (1999) Can one hear the
shape of an electrode? The active zone in Laplacian transfer, II.
Theory. {\it European Journal of Physics,  B, \bf 9}, 754--763.

\item{[Fu]} M.~Fukushima (1967), A construction of reflecting barrier Brownian
motions for bounded domains. {\it Osaka J. Math.}, {\bf 4},
183-215.

\item{[FOT]} M.~Fukushima, Y.~Oshima and M.~Takeda (1994),
 {\it Dirichlet Forms and Symmetric Markov Processes}.
 de Gruyter, Berlin.

\item{[FT]} M. Fukushima and M. Tomisaki (1996), Construction and
decomposition of reflecting diffusions on Lipschitz domains with
H\"older cusps. {\it Probab. Theory Rel. Fields \bf 106},
521--557.

\item{[GFS]} D.S.~Grebenkov, M.~Filoche and B.~Sapoval (2003)
Spectral properties of the Brownian self-transport operator {\it
Europhys. Journ. B. \bf 36}, 221--231.

\item{[GA]} K.~Gustafson and T.~Abe (1998), The third boundary
condition---was it Robin's? {\it Math. Intelligencer \bf 20},
63--71.

\item{[K]} C. E. Kenig (1994), {\it Harmonic Analysis Techniques for Second
Order Elliptic Boundary Value Problems}. CBMS 83,
Math. Assoc. Amer..

\item{[LSW]} W. Littman, G.Stampacchia, and H.F. Weinberger (1963),
Regular points for elliptic equations with discontinuous coefficients.
{\it Ann. Scuola  Norm. Sup. Pisa \bf 17}, 43--77.

\item{[Mv]} B.~Maisonneuve (1975), Exit systems, {\it Ann.
Probability \bf 3}, 399--411.

\item{[MS]} Z.-M. Ma and R. Song (1990), Probabilistic methods in
Schr\"odinger equations. {\it Seminar on Stochastic Processes,
1989}. 135--164, {\it Progr. Probab., \bf 18}, Birkh\"auser,
Boston.

\item{[Maz]} V.~G.~Maz'ja (1985), {\it Sobolev Spaces}. Springer Series in
Soviet Mathematics. Springer-Verlag, Berlin.

\item{[Pa]} V.G.~Papanicolaou (1990), The probabilistic solution of
the third boundary value problem for second order elliptic
equations. {\it Probab. Theory Rel. Fields \bf 87}, 27--77.

\item{[Po]} Ch. Pommerenke (1992), {\it Boundary Behaviour of Conformal
Maps}, Springer, Berlin.

\item{[Sa]} B.~Sapoval (1996), Transport across irregular interfaces:
fractal electrodes, membranes and  catalysts, in {\it Fractals and
disordered systems, 2nd ed.}, A.~Bunde and S.~Havlin, Eds.
(Springer-Verlag) 232--261.

\item{[Sh]} M.~Sharpe (1988), \it General Theory of Markov
Processes, \rm Academic Press, Boston.

\item{[Sm1]} R. Smits (2005), Eigenvalue estimates and critical temperature
in zero fields for enhanced surface superconductivity, preprint.

\item{[Sm2]} R. Smits (2005), Upper and lower bounds for the
principal eigenvalue of the generalized Robin problem, preprint.

\vskip1truein

\noindent R.B.: Department of Mathematics, University of
Connecticut, Storrs, CT 06269-3009 {\tt bass@math.uconn.edu}

\bigskip
\noindent K.B. and Z.-Q. C.: Department of Mathematics, Box
354350, University of Washington, Seattle, WA 98115-4350, USA
\hfill\break{\tt burdzy@math.washington.edu,
zchen@math.washington.edu}

\bye